\pdfoutput=1
\RequirePackage{ifpdf}
\ifpdf 
\documentclass[pdftex]{sigma}
\else
\documentclass{sigma}
\fi

\usepackage[all]{xy}

\newtheorem{Theorem}{Theorem}[section]
\newtheorem{Corollary}[Theorem]{Corollary}
\newtheorem{Proposition}[Theorem]{Proposition}
 { \theoremstyle{definition}
\newtheorem{Definition}[Theorem]{Definition}
\newtheorem{Example}[Theorem]{Example} }

\numberwithin{equation}{section}

\begin{document}


\newcommand{\arXivNumber}{1308.3819}

\renewcommand{\PaperNumber}{084}

\FirstPageHeading

\ShortArticleName{Fast Basins and Branched Fractal Manifolds of Attractors of Iterated Function Systems}

\ArticleName{Fast Basins and Branched Fractal Manifolds\\ of Attractors of Iterated Function Systems}

\Author{Michael F.~{BARNSLEY}~$^\dag$ and Andrew {VINCE}~$^\ddag$}

\AuthorNameForHeading{M.F.~Barnsley and A.~Vince}

\Address{$^\dag$~Mathematical Sciences Institute, Australian National University, Australia}
\EmailD{\href{mailto:michael.barnsley@anu.edu.au}{michael.barnsley@anu.edu.au}}
\URLaddressD{\url{http://www.superfractals.com}}

\Address{$^\ddag$~Department of Mathematics, Univesity of Florida, USA}
\EmailD{\href{mailto:avince@ufl.edu}{avince@ufl.edu}}

\ArticleDates{Received June 23, 2015, in f\/inal form October 13, 2015; Published online October 16, 2015}

\Abstract {The fast basin of an attractor of an iterated function system (IFS) is the
set of points in the domain of the IFS whose orbits under the associated
semigroup intersect the attractor. Fast basins can have non-integer
dimension and comprise a class of deterministic fractal sets. The
relationship between the basin and the fast basin of a point-f\/ibred
attractor is analyzed. To better understand the topology and geometry of fast basins, and because
of analogies with analytic continuation, branched fractal manifolds are
introduced. A branched fractal manifold is a metric space constructed from
the extended code space of a point-f\/ibred attractor, by identifying some
addresses. Typically, a branched fractal manifold is a union of a
nondenumerable collection of nonhomeomorphic objects, isometric copies of
generalized fractal blowups of the attractor.}


\Classification{05B45; 37B50; 52B50}

\section{Introduction}\label{sec:intro}

This paper on iterated function systems concerns related concepts: the fast
basin, fractal continuation, an extended code space, and the branched
fractal manifold.

The fast basin of an attractor of an iterated function system (IFS) is the
set of points, each of which possesses a chaos game orbit \cite{chaos} that
reaches the attractor in f\/initely many steps. More precisely, let $A$ be an
attractor of an IFS $\mathcal{F}$, here def\/ined to be a f\/inite set of
homeomorphisms $f_{n}\colon \mathbb{X}\rightarrow \mathbb{X}$, $n=1,2,\dots ,N$,
mapping a complete metric space~$\mathbb{X}$ onto itself. If $f,g\in
\mathcal{F}$, then $f\circ g$ denotes the composition of $f$ with $g$. The
fast basin of $A$ is the set of points $x\in \mathbb{X}$ such that there
exists a f\/inite sequence of digits $\{\iota _{n}\}_{n=1}^{k}\subset
\{1,2,\dots ,N\}^{k}$ with $f_{\iota _{1}}\circ f_{\iota _{2}}\circ \cdots
\circ f_{\iota _{k}}(x)\in A$. Equivalently, the fast basin is the
complement of the set of points, all of whose orbits do not intersect the
attractor. If the maps that comprise the IFS are contractive similitudes on~$\mathbb{R}^{n}$, then the fast basin is the union of fractal blowups, as def\/ined by Strichartz \cite{strichartz0}, of $A$ with respect to the IFS. In general, a fast basin is a union of fractal continuations, generalizing the main idea in~\cite{continuations}, of $A$ with respect to the IFS.

Some examples of fast basins, corresponding to geometrically simple IFSs,
are presented in Section~\ref{examplesec}, illustrating that they can have interesting
geometrical structure and that they are not easy to understand. Fast basins
of af\/f\/ine IFSs may have non-integer Hausdorf\/f dimension, as proved in~\cite{mfbkl}, and
thus comprise a class of deterministic fractals which may be suitable for
geometrical modelling of real-world objects, extending the reach of ideas initiated by Mandelbrot; see for example the many modelling papers that cite~\cite{mandelbrot} and use IFS theory. But our interest here is in mathematical structure of fast basins, which is relevant to fractal transformations~\cite{monthly}, tiling theory~\cite{tiling}, non-commutative geometry~\cite{ionescu}, and to a generalization of analytic continuation (see below).

\looseness=-1
A branched fractal manifold is a certain metric space built from the code
space of a point-f\/ibred attractor of an IFS. Its projection onto~$\mathbb{X}$
is the fast basin of the attractor. An initial motivation for def\/ining and
analyzing branched fractal manifolds was a desire to understand the
topological structures of fast basins, at least in the case of simple IFSs
consisting of a few af\/f\/ine maps.

A strong motivation for interest in fast basins and branched fractal
manifolds is a new view of analytic continuation, introduced in~\cite{continuations}. According to this view, fractal continuation provides an
extension
of the notion of analytic continuation. As a simple example, consider the
IFS def\/ined by the pair of M\"{o}bius transformations,
\begin{gather*}
\begin{pmatrix}
-31+4i & 8+22i \\
2+11i & 2-4i%
\end{pmatrix}\qquad
\text{and}\qquad
\begin{pmatrix}
-25-13i & -17+14i \\
-11+7i & -4+13i%
\end{pmatrix}%
\end{gather*}
acting on the Riemann sphere. This IFS has a unique attractor, an arc of a
circle (of radius one half, centered at~$3i/2$). The fast basin is the whole
circle, i.e., the analytic continuation of the arc. This illustrates a more
general situation. In~\cite{continuations} it is proved, under conditions
that are not stated here, that the graph of a real analytic function over an
interval, say $[0,1]$, is the attractor of an IFS whose maps are analytic on
a neighborhood of the graph, and generalized blowups of the graph, generated
using the IFS, coincide with the graph of any analytic continuation of the
original function, where both are def\/ined. (In~\cite{continuations} and here,
we refer to generalized fractal blowups of the graphs of fractal
interpolation functions as (fractal) continuations.) In~\cite{continuations}
it is further established that all except two fractal continuations of the
graphs of any af\/f\/ine fractal interpolation function $g\colon [0,1]\rightarrow
\mathbb{R}$ are def\/ined for all real~$x$. Although a fractal continuation
coincides with the real analytic continuation of~$g$ when~$g$ is real
analytic on a neighborhood of $[0,1]$, in general the graphs of dif\/ferent
continuations of an af\/f\/ine fractal interpolation function do not coincide,
except over $[0,1]$. Instead, their union, which is an example of a fast
basin, may be interpreted as an analog or generalization of analytic
continuation applied to a fractal function~$g$. In this interpretation,
dif\/ferent continuations of $g$ belong to dif\/ferent branches, thus the need
to understand how these branches f\/it together and the motivation for the
fast basin and branched fractal manifold.

The organization of the paper is as follows. In Section~\ref{sec:IFS} we
def\/ine point-f\/ibred IFSs, their attractors, duals, basins, fast basins and
fractal continuations. A~fractal continuation is a~ge\-ne\-ralization of a~fractal blowup as def\/ined by Strichartz~\cite{strichartz0}. Proposition~\ref{prop:FB} states that the fast basin of an attractor is the union of its
fractal continuations. Section~\ref{examplesec} provides examples of fast
basins.

The relationship between the basin and the fast basin of an attractor is the
subject of Section~\ref{sec:BFB}. Theorem~\ref{reversethm} shows that the
image of a disjunctive point in code space under the coding map (see below)
belongs to the boundary of an attractor if and only if the interior of the
attractor is empty. Theorem~\ref{interiorthm} relates addresses of points in
attractors to whether or not associated continuations contain the basin,
using the concept of reversible addresses. As a consequence, Corollary~\ref{cor:basins} states that the fast
basin contains the basin if and only if the attractor has nonempty interior.

In Section~\ref{symbtheorysec}, two symbolic IFSs $Z$ and $\widehat{Z}$ are introduced. The
f\/irst, $Z=\{\mathbb{I}_{+};\,\sigma_{n},\,n=1,2,\dots ,N\}$, involves shift maps on
$\mathbb{I}_{+}=\{1,2,\dots ,N\}^{\infty}$ and is well known, when~$\mathcal{F}$
is contractive, to be related to the attractor $A$ by the following
commutative diagram
\begin{gather*}  
\xymatrix{
  \mathbb{I}_{+} \ar[r]^{\sigma_{n}}    \ar[d]_{\pi} & \mathbb{I}_{+}  \ar[d]^{\pi} \\
        A  \ar[r]_{f_{n}}  &   A   }
\end{gather*}
This involves the coding map or addressing function, a continuous surjection
$\pi\colon \mathbb{I}_{+}\rightarrow A$.
However, $Z$ is not an IFS, as
def\/ined in this paper, as the maps of $Z$ are not invertible and its fast
basin is undef\/ined. This is remedied by the introduction of $\widehat{Z}=\{\mathbb{I};\,\sigma_{n},\,n=1,2,\dots ,N\}$ whose maps are invertible. Here $\mathbb{I}$ is a certain shift invariant subspace of $\mathbb{I}_{0} =\{-N,
\dots,-1,1,2, \dots ,N\}^{\infty}$. The attractor of $\widehat{Z}$ is $\mathbb{I}_{+}$ and its fast basin, $\widehat{\mathbb{I}}\subset\mathbb{I}$,
provides a symbolic representation for any attractor of $\mathcal{F}$ and
its fast basin. It is shown that the coding map $\pi\colon \mathbb{I}_{+}\rightarrow A$ can be extended to various subspaces of~$\mathbb{I}$,
including~$\widehat {\mathbb{I}}$, and that the action of shift maps on
these subspaces is (semi-)conjugate to the action of the~$f_{n}$s (and/or
their inverses) on corresponding subspaces of~$\mathbb{X}$, including the
fast basin of~$A$. Dual relationships also hold in certain cases. These
results are summarized by the commutative diagrams in Theorem~\ref{extenthm}
and Corollaries~\ref{theorem1} and~\ref{dualcasethm}.

The branched fractal manifold $\mathbb{L}$ associated with a point-f\/ibred
attractor of an invertible IFS is def\/ined in Section~\ref{leafsec}. It is a
metric space constructed from the symbolic fast basin $\widehat {\mathbb{I}}$
by identifying some addresses. It is a disjoint union of certain leaves,
each leaf being homeomorphic to one of at most~$N+1$ subsets of~$A$, but~$\mathbb{L}$ itself may not be locally compact. Also, $\mathbb{L}$ contains
isometric copies, called sheets, of all the fractal continuations of the
original attractor. The branched fractal manifold contains, in a concise
way, information on the global addressing discussed in Section~\ref{symbtheorysec}. A description of~$\mathbb{L}$, providing its main
properties, is the content of Theorems~\ref{metricthm} and~\ref{sevenpttwothm}. Examples are provided.

Relationships between this paper and related work in the literature are
discussed in Section~\ref{sec:literature}.

\section{IFSs attractors, their basins, fast basins, and
continuations}\label{sec:IFS}

Let $\mathbb{N}=\{1,2,3,\dots \}$ and $\mathbb{N}_{0}=\{0,1,2,\dots \}$.
Throughout this paper an \textit{iterated function system $($IFS$)$} is a~complete metric space $\mathbb{X}$ together with a f\/inite set of
homeomorphisms $f_n\colon \mathbb{X}\rightarrow \mathbb{X}$, $n = 1,2, \dots, N$,
and is denoted by
\begin{gather*}
\mathcal{F}=\mathcal{F}_{\mathbb{X}}=\{\mathbb{X};\, f_{1},f_{2},\dots ,f_{N}\}.
\end{gather*}
We use the same symbol $\mathcal{F}$ for the IFS and for the set of
functions $\{f_{1},f_{2},\dots ,f_{N}\}$.

Let $H=H_{\mathbb{X}}$ be the collection of nonempty compact subsets of $\mathbb{X}$ and def\/ine $F\colon H\rightarrow H$ by
\begin{gather*}
F(C)=\bigcup\limits_{f\in\mathcal{F}}f(C)
\end{gather*}
for all $C\in H$, where $f(C)=\{f(x)\colon x\in C\}$. We extend naturally $F$ to a map on the collection of all
subsets of $\mathbb{X}$. For $S\subset\mathbb{X}$, def\/ine $F^{0}(S)=S$ and
let~$F^{k}(S)$ denote the $k$-fold composition of~$F$ applied to $S$,
namely, the union of $\{\iota_{n}\}_{n=1}^{k}\subset\{1,2,\dots ,N\}^k$ with $%
f_{\iota_{1}}\circ f_{\iota_{2}}\circ\cdots\circ f_{\iota_{k}}(S)$ over all
f\/inite words $\iota_{1}\iota_{2}\cdots \iota_{k}$ of length~$k$.

Let $d=d_{\mathbb{X}}$ be the metric on $\mathbb{X}$, and let $d_{H}=d_{H_{\mathbb{X}}}$ be the corresponding Hausdorf\/f metric on~$H_{{\mathbb{X}}}$.
Throughout, the topology on~$H_{\mathbb{X}}$ is the one induced by~$d_{H_{\mathbb{X}}}$. It is well known that $(H,d_{H})$ is a complete metric space
because~$(\mathbb{X},d)$ is complete, and that if~$(\mathbb{X},d)$ is
compact then~$(H,d_{H})$ is compact.

\begin{Definition}
\label{def:attractor} An \textit{attractor} of the IFS $\mathcal{F}$ is a
set $A\in H$ such that
\begin{enumerate}\itemsep=0pt
\item[1)] $F(A)=A$, and

\item[2)] there is an open set $U\subset\mathbb{X}$ such that $A\subset U$ and $\lim\limits_{k\rightarrow\infty}F^{k}(C)=A$, for all $S\in H$ with $C \subset U$,
where the limit is w.r.t.\ the metric $d_{H}$.
\end{enumerate}
\end{Definition}

The union of all open sets $U$, such that statement~2 of Def\/inition~\ref{def:attractor} is true, is called the \textit{basin }of the attractor~$A$
(w.r.t.~$\mathcal{F}$). If $B=B(A)$ denotes the basin of~$A$, then it can be
proved that statement~2 of Def\/inition \ref{def:attractor} holds with~$U$
replaced by~$B$. That is, the basin of the attrac\-tor~$A$ is the largest open
set~$U$ such that statement~2 of Def\/inition~\ref{def:attractor} holds. An
IFS may not possess an attractor, or it may possess multiple attractors.
Examples and further discussion can be found in~\cite{BV8} and in references
therein.

The IFS $\mathcal{F}$ is said to be \textit{contractive} if there is a $\lambda \in [0,1)$ such that
\begin{gather*}
d_{\mathbb{X}}(f(x),f(y))\leq\lambda d_{\mathbb{X}}(x,y)
\end{gather*}
for all $f\in \mathcal{F}$ and all $x,y\in\mathbb{X}$. A basic result of
\cite{Hutchinson} is that a contractive IFS possesses a unique attractor~$A$, with $B(A)=\mathbb{X}$.

Note that, in this paper, each $f_n\colon \mathbb{X}\rightarrow \mathbb{X}$ is a homeomorphism and so has a continuous inverse $f_n^{-1}\colon \mathbb{X}\rightarrow \mathbb{X}$ which is also a homeomorphism. This allows us to def\/ine an ``inverse'' of the IFS~$\mathcal{F}$ as follows.

\begin{Definition}
The \textit{dual IFS} is
\begin{gather*}
\mathcal{F}^{\ast} :=\big\{\mathbb{X};\, f_{n}^{-1},\, n=1,2,\dots ,N\big\}.
\end{gather*}
\end{Definition}

If $A$ is an attractor of the IFS $\mathcal{F}$, then the set $A^{\ast}:=
\mathbb{X}\backslash B(A)$ is called the \textit{dual repeller} of~$A$ (w.r.t.~$\mathcal{F}$). The dual repeller~$A^{\ast}$, under certain
conditions \cite{BVaustMS}, is an attractor of the dual IFS~$\mathcal{F}^{\ast}$. This can occur, for example, when $\mathcal{F}_{\widehat{\mathbb{C}
}}$ is a M\"{o}bius IFS on the Riemann sphere~$\widehat{\mathbb{C}}$~\cite{Vince}.

For an inf\/inite word $\iota = \iota_{1} \iota_{2} \iota_{3} \cdots$, $\iota_i \in \{1,2, \dots ,N\}$, let
\begin{gather}  \label{eq:composition}
  \iota|k  = \iota_{1} \iota_{2} \iota_{3} \cdots \iota_k, \qquad
f_{\iota|k}  = f_{\iota_{1}}\circ f_{\iota_{2}}\circ \cdots \circ
f_{\iota_k},
\end{gather}
and $f_{\varnothing}=f_{\iota|0}={\text{id}}_{\mathbb{X}}$, the identity map
on~$\mathbb{X}$.

\begin{Definition}
\label{def:PF} An attractor $A$ of $\mathcal{F}$ is \textit{point-fibred}
(w.r.t.~$\mathcal{F}$) if
\begin{gather*}
\lim_{k\rightarrow\infty}f_{\iota|k}(C) \subset A,
\end{gather*}
is a singleton subset of $\mathbb{X}$, for all $\iota \in \{1,2, \dots
,N\}^{\infty}$ independent of $C\subset B(A)$ with $C\in H$, where
convergence is with respect to the Hausdorf\/f metric.
\end{Definition}

All attractors in this paper are assumed to be point-f\/ibred. It follows from
results in \cite{Hutchinson} that if the IFS $\mathcal{F}$ is contractive,
then its attractor is point-f\/ibred. The point-f\/ibred property allows for an
addressing scheme for points of the attractor (Def\/inition~\ref{def:CM}). The
following notation will be used throughout this paper. Let $N\in\mathbb{N}$
and let
\begin{gather}
 I_{+}  =\{+1,+2, \dots ,+N\},\qquad \mathbb{I}_{+}  =I_{+}^{\infty}, \nonumber\\
 I_{-}   =\{-1,-2, \dots ,-N\},\qquad \mathbb{I}_{-}  =I_{-}^{\infty}, \nonumber\\
 I   =I_{+}\cup I_{-},   \qquad \hspace{20mm}\mathbb{I}_{0}= I^{\infty},\nonumber\\
  \label{eq:notation}
\mathbb{I=\{\iota\in I}_{0}\colon \iota_{i}\neq-\iota_{i+1}\text{ for all } i\in \mathbb{N\}}.
\end{gather}
Let $d_{\mathbb{I}}$ be the metric on $\mathbb{I}_{0}$ def\/ined by
\begin{gather*}
d_{\mathbb{I}}(\iota,\omega)=
\begin{cases}
\max\big\{2^{-k}\colon k\in\mathbb{N}, \iota_{k}\neq\omega_{k}\big\} & \text{if} \ \
\iota \neq\omega, \\
0 & \text{if} \ \ \iota =\omega,
\end{cases}
\end{gather*}
for all $\iota,\omega\in\mathbb{I}_{0}$. The metric $d_{\mathbb{I}}$ induces
the product topology on $\mathbb{I}_{0}$, and the metric spaces $(\mathbb{I}, d_{\mathbb{I}})$,
 $(\mathbb{I}_{+}, d_{\mathbb{I}}),(\mathbb{I}_{-},d_{\mathbb{I}})$ are compact subspaces of the compact metric space $(\mathbb{I}_{0},d_{\mathbb{I}})$.

The space $(\mathbb{I}_{+},d_{\mathbb{I}})$ is sometimes called the \textit{code space} or \textit{address space} for an attractor of~$\mathcal{F}$
because it provides addresses for the points of~$A$, as given in the
following def\/inition.

\begin{Definition}\label{def:CM}
Let $A$ be a point-f\/ibred attractor of an IFS $\mathcal{F}$.
According to Def\/inition~\ref{def:PF}, there is a map $\pi\colon \mathbb{I}_{+}\rightarrow A\subset\mathbb{X}$, called the \textit{coding map}, that is
well def\/ined by
\begin{gather}
\pi(\iota)=\lim_{k\rightarrow\infty}f_{\iota|k}(b),  \label{limiteqn}
\end{gather}
for all $\iota\in\mathbb{I}_{+}$, and the limit is independent of $b\in B(A)
$. It is readily proved that $\pi$ is continuous and that $\pi(\mathbb{I}_{+})=A$. The set-valued inverse $\pi^{-1}(a)$ comprises the \textit{set of
addresses} of the point $a\in A$ (w.r.t.~$\mathcal{F}$).
\end{Definition}

\begin{Definition}
The \textit{fast basin} $\widehat{B}=\widehat{B}(A)$ of an attractor $A$ of
an IFS $\mathcal{F}$ is
\begin{gather*}
\widehat{B}=\big\{x\in\mathbb{X} \colon F^{K}(\{x\})\cap A\neq\varnothing \
\text{for some} \  K \in \mathbb{N }\big\}.
\end{gather*}
\end{Definition}

The relationship between the basin $B$ and the fast basin $\widehat{B}$ is
complicated because neither $B\subset$ $\widehat{B}$ nor $\widehat{B}\subset
B$ holds in general. This topic is discussed in Section~\ref{sec:BFB}.

Fractal continuations were introduced in~\cite{continuations} in the context
of fractal interpolation, and in~\cite{BV8} in the context of fractal
tiling. (See also the last part of Section~\ref{sec:intro} for a discussion
of related notions.) These situations are generalized in Def\/inition~\ref{def:continuation}. Extend the notation of equation~\eqref{eq:composition}
as follows. For $\iota \in \mathbb{I}_{0}$ def\/ine $-\iota$ by
\begin{gather*}
(-\iota)_k := -\iota_k
\end{gather*}
for all $k \in \mathbb{N}$ and
\begin{gather*}
f_{-j} = (f_j)^{-1}
\end{gather*}
for $j = 1,2,\dots, N$.

\begin{Definition}\label{def:continuation} Let $A$ be an attractor of an IFS $\mathcal{F}=\{\mathbb{X} ;\,f_{n},\,n\in I_{+}\}$.
For $\theta\in\mathbb{I}_{+}$, the \textit{continuation} ${B}_{\theta }$ of $A$ (w.r.t.~$\mathcal{F}$) is def\/ined to be
\begin{gather*}
B_{\theta} =\bigcup\limits_{k\in\mathbb{N}} f_{-\theta|k}(A) = A \cup
\bigcup\limits_{k\in\mathbb{N}} f_{\theta_{1}} ^{-1}\circ
f_{\theta_{2}}^{-1}\circ\dots\circ f_{\theta_{k}}^{-1}(A).
\end{gather*}
Note that this is a nested union in the sense that $f_{-\theta|k}(A) \subset
f_{-\theta|k+1}(A)$.
\end{Definition}

The set ${B}_{\theta}$ is referred to as a~\textit{continuation} of~$A$ or a~\textit{fractal continuation} of~$A$. The family $\{{B}_{\theta }\colon \theta\in \mathbb{I}_{+}\}$ is referred to as the set of continuations of~$A$. We also
write ${B}_{\theta |k}=f_{-\theta|k}(A)$, which is referred to as a \textit{finite continuation of~$A$}. Clearly
\begin{gather*}
{B}_{\theta}=\bigcup_{k\in\mathbb{N}_{0}}{B}_{\theta|k}.
\end{gather*}
The paper \cite{continuations} concerns the special case of fractal
continuations of fractal functions. In that case it is proved, under special
conditions, the most important of which is that the IFS consists of analytic
functions, that the set of continuations is uniquely determined by the
attractor, independent of the analytic IFS used to generate the attractor.

The following result, stating that the fast basin is the union of all the
continuations, is readily verif\/ied.

\begin{Proposition}\label{prop:FB} Let $A$ be an attractor of an IFS $\mathcal{F} =\{\mathbb{X};\,f_{n},\,n\in I_{+}\}$. If $\widehat{B}$ is the fast basin of $A$ and $\{{B}_{\theta}\colon \theta\in\mathbb{I}_{+}\}$ is the set of continuations of~$A$,
then
\begin{gather*}
\widehat{B}=\bigcup\limits_{\theta\in {\mathbb{I}}_{+}}{B}_{\theta}.
\end{gather*}
\end{Proposition}

\section{Fast basin examples}\label{examplesec}

\begin{Example}[Cantor set] \label{ex:CS} Routine calculation shows that the fast
basin of the standard Cantor set $\mathcal{C}\subset\mathbb{R}$, w.r.t.\ the
IFS $\left\{ \mathbb{R} ;\,x/3,\,x/3+2/3\right\}$, is
\begin{gather*}
\widehat{B}=\bigcup_{k\in\mathbb{Z}}(\mathcal{C+}2k).
\end{gather*}
We use the notation for Minkowski sum, $(\mathcal{C+}2k):=\{x+2k\colon x\in \mathcal{C}\}$.
\end{Example}

\begin{Example}[Sierpinski gasket] The fast basin of a Sierpinski triangle $\triangle\subset \mathbb{R}^{2}$, with vertices at $v_1,v_2,v_3 \mathbb{\in R}^{2}$, w.r.t.\ the IFS $\left\{ \mathbb{R}^{2};\,(x+v_1)/2,\,(x+v_2)/2,\,(x+v_3)/2\right\} $, is $\widehat {B}=\cup_{t\in\mathcal{G}}(\triangle+t)$ where $\mathcal{G}$ is the group
generated by the set of translations by $\{v_1-v_2,v_2-v_3,v_3-v_1\}$. This may be proved by induction, using $A=\cup_i f_i(A)$ and $f_{v_i}^{-1}\circ f_{v_j}$ is translation by $v_j-v_i$. Fig.~\ref{manifold} shows~$\widehat{B}$ when $v_1=(0,0)$, $v_2=(1,0)$, and $v_3=(0,1)$.

\begin{figure}[t]
\centering
\fbox{\includegraphics[width=2.919in]{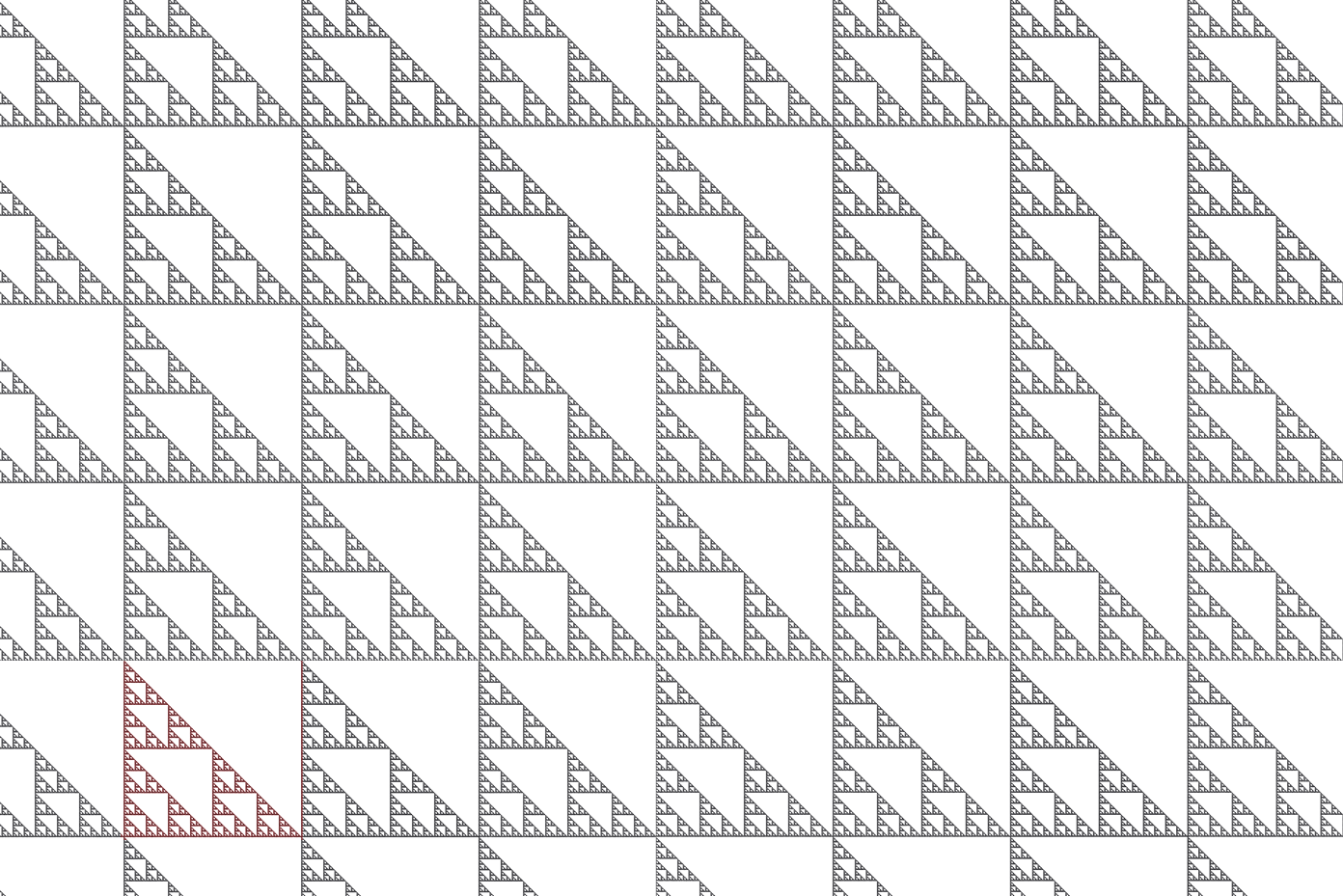}}
\caption{This illustrates part of the fast basin of the attractor of the
iterated function system $\mathcal{F}=\{\mathbb{R}^{2};\,(x/2,y/2),\,((x+1)/2,y/2),\,(x/2,(y+1)/2)\}$. The attractor is rendered in
red at the lower left.}
\label{manifold}
\end{figure}
\end{Example}

\begin{Example}[other af\/f\/ine IFSs]
\label{ex:affine} Fig.~\ref{kockcurve3lr}
illustrates part of the fast-basin of a contractive af\/f\/ine IFS $\{ {\mathbb{R}} ^2 ; \, f_1,\, f_2 \}$, where
\begin{gather*}
 f_1(x,y) = \left (\frac{x}{2}+\frac{y}{2\sqrt{3}} - 1,
\frac {x}{2\sqrt{3}}-\frac{y}{2} \right ), \qquad
 f_2(x,y) = \left(\frac{x}{2}-\frac{y}{2\sqrt{3}} +1,  -\frac {x}{2\sqrt{3}}-\frac{y}{2}
\right ).
\end{gather*}
The attractor is the segment of the Koch snowf\/lake curve near the center of
the f\/igure. The fast basin is the union of the boundaries of the tiles of a
tiling of the plane by Koch snowf\/lakes and other related tiles.
\begin{figure}[t]
\centering
\fbox{\includegraphics[natheight=7.111400in,
natwidth=7.111400in,
height=2.0358in,
width=2.0358in]{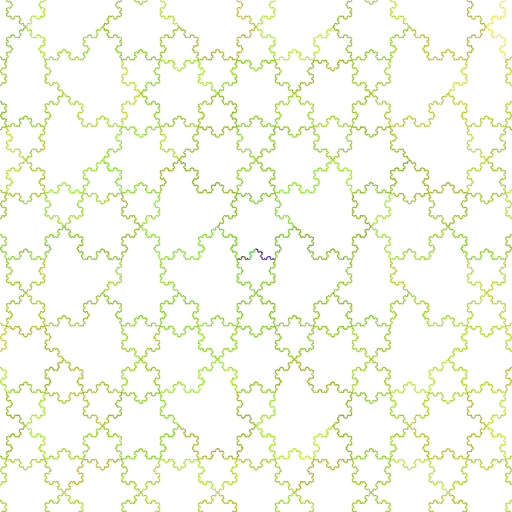}}
\caption{See Example~\protect\ref{ex:affine}. This shows the fast basin for
the Koch curve using an IFS of similitudes. The Koch curve is located near
the center, in various colours.}\label{kockcurve3lr}
\end{figure}

Fig.~\ref{parab4fract} illustrates part of the fast basin for the
contractive af\/f\/ine IFS $\{{\mathbb{R}}^2 ; \, f_1,\, f_2\}$, where
\begin{gather*}
  f_1(x,y)  = \left ( \frac{x}{2}+\frac{1}{2},\frac{x}{2}+\frac{2y}{5}+\frac {1}{4} \right ), \qquad
f_2(x,y)  = \left (\frac{x}{2} -\frac{1}{2},  -\frac{x}{2}+\frac{2y}{5}+\frac {1}{4} \right ),
\end{gather*}
whose attractor, at the center of the image, is the graph of a fractal
interpolation function. At each branch point, there is a countable inf\/inity
of distinct branches.

Both f\/igures in this example were obtained by direct computation.

\begin{figure}[t]
\centering
\fbox{\includegraphics[natheight=7.111400in,
natwidth=7.111400in,
height=2.0358in,
width=2.0358in]{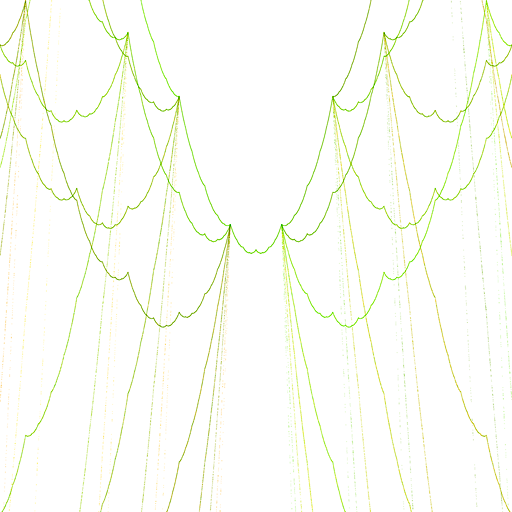}}
\caption{See Example~\protect\ref{ex:affine}.}\label{parab4fract}
\end{figure}
\end{Example}

\begin{Example}[topological dimension of the fast basin] We use the following
def\/inition: an \textit{orbit of a point} $x\in\mathbb{X}$, under
the IFS $\mathcal{F} =\{\mathbb{X};\, f_{1}, \dots,f_{N}\}$, is
$\left\{ f_{\iota_{1}}\circ f_{\iota_{2}} \circ \cdots \circ
f_{\iota_{k}}(x)\right\} _{k=1}^{\infty}\subset\mathbb{X}$ for some $\iota\in\left\{ 1,2,\dots ,N\right\} ^{\infty}$. If $\mathcal{F}$ is an IFS on ${\mathbb{R}}^M$ consisting of af\/f\/ine functions and whose attractor $A$ is
contained in a proper linear subspace~$V$ of $\mathbb{R}^{M}$, then it is
easy to show that the fast basin of $A$ is also contained in~$V$. Moreover,
any orbit of any point not in $V$ has empty intersection with~$A$. It
follows that the topological dimension of the fast basin of $A$ is strictly
less than that of the underlying space. For instance, if $M=N=2$, and
\begin{gather*}
f_{1}(x,y)=(x/2,y/2),\qquad  f_{2}(x,y)=((x+1)/2,y/2)
\end{gather*}
for all $(x,y)\in\mathbb{R}^{2}$, then $A=[0,1] \times \{0\} \subset\mathbb{R}\times \{0\}$, $B= \mathbb{R}^{2}$, and $\widehat{B}\subset\mathbb{R}\times
\{0\}$. Any orbit of any point in $\mathbb{R}^{2}{\setminus} \{(x,0)\colon x
\in {\mathbb{R}} \}$ does not meet the attractor (although the closure of
the orbit does meet the attractor) while, for any point in~$\mathbb{R}$,
there exists an orbit that reaches $A$ in f\/initely many steps, i.e., the
orbit intersects the attractor.
\end{Example}

\begin{Example}
An example, illustrating a relationship between a fast basin and analytic
continuation, is provided by the contractive IFS
\begin{gather*}
\mathcal{F} =\big\{\mathbb{C}^2; \, f_{+1+i},\,f_{-1+i},\,f_{+1-i},\,f_{-1-i}\big\},
\end{gather*}
where
\begin{gather*}
f_{\pm1-i}
\begin{pmatrix}
z \\
w
\end{pmatrix}
 =
\begin{pmatrix}
\frac{1}{2} & 0 \\
\frac{\pm1-i}{2} & \frac{1}{4}
\end{pmatrix}
\begin{pmatrix}
z \\
w
\end{pmatrix}
+
\begin{pmatrix}
\frac{\pm1-i}{2} \\
\frac{\mp i}{2}
\end{pmatrix}
, \\
f_{\pm1+i}
\begin{pmatrix}
z \\
w
\end{pmatrix}
 =
\begin{pmatrix}
\frac{1}{2} & 0 \\
\frac{\pm1+i}{2} & \frac{1}{4}
\end{pmatrix}%
\begin{pmatrix}
z \\
w
\end{pmatrix}
+
\begin{pmatrix}
\frac{\pm1+i}{2} \\
\frac{\pm i}{2}
\end{pmatrix}.
\end{gather*}
By following arguments similar to ones in \cite{continuations}, we f\/ind that the unique attractor $A$ of $\mathcal{F}$ is the graph of $z\mapsto z^{2}$
over the square $-1\leq\operatorname{Re}z,\operatorname{Im}z\leq+1$. We also f\/ind that the fast basin~$\widehat{B}$ of $A$ (w.r.t.~$\mathcal{F}$) is the manifold
\begin{gather*}
\big\{ \big(z,z^{2}\big) \colon z\in\mathbb{C}\big\} .
\end{gather*}
\end{Example}

\begin{Example}\label{ex:schottky1} An example that is related to a Schottky group is
provided by the loxodromic M\"{o}bius IFS
\begin{gather*}
\mathcal{F} =\left \{\widehat{\mathbb{C}};\, f_{1}(z)=\frac{az-i(2+\sqrt{3})
}{z-i(2+\sqrt{3})+a-1},\, f_{2}(z)=\frac{az+i(2+\sqrt{3})}{-z-i(2+\sqrt{3}
)+a-1}\right \}
\end{gather*}
where $a=\frac{-i(2+\sqrt{3})C+1}{1-C}$, $C\in\mathbb{C}$ with $\left\vert
C\right\vert \notin\{0,1\}$, $f_{1}(1)=1$ and $f_{2}(-1)=-1$. We f\/ind
\begin{gather*}
f_{1}^{\prime}(1)=f_{2}^{^{\prime}}(-1)=C,
\end{gather*}
so if $\left\vert C\right\vert $ is suf\/f\/iciently small, then $\mathcal{F}$ has a totally
disconnected attractor $A$, located within two small circles, one centered
at $+1$ and one centered at $-1$. This is because the center of each circle is an attractive f\/ixed point of one of the maps, and the derivative of each map, at its attractive f\/ixed point, has magnitude $\left\vert C\right\vert $.  Similarly, in this case, the dual repeller~$A^{\ast}$
is also totally disconnected and is located within two small circles, one
centered at $(2+\sqrt{3})i$ and the other at $-(2+\sqrt{3})i$. The maps $f_{1}$, $f_{2}$ are the generators of a Schottky group whose limit set~$\Lambda$
is also totally disconnected, and we have $A\subset\Lambda$, $A^{\ast}\subset\Lambda$, and $\widehat{B}\subset\Lambda$. It follows that $\widehat{B}$ is totally disconnected.
\end{Example}

\section{When does the fast basin contain the basin?}\label{sec:BFB}

Consider the following examples showing that, in general, there is no
containment relationship between the basin and the fast basin. Since the IFS
\begin{gather*}
\left \{ {\mathbb{R}}; \, \frac12  x, \,\frac12   x + \frac12\right \}
\end{gather*}
is contractive, the basin $B = {\mathbb{R}}$. It is not hard to show, for
example via Proposition~\ref{prop:FB}, that the fast basin $\widehat B = {\mathbb{R}}$. Therefore, $\widehat B = B$ for this IFS. Likewise, the IFS
\begin{gather*}
\left \{ \mathbb{R}; \, \frac13   x, \, \frac13   x+\frac{2}{3}\right \}
\end{gather*}
is contractive, but, as in Example~\eqref{ex:CS}, the fast basin is not all
of ${\mathbb{R}}$. Therefore $\widehat B \subsetneq B$ for this IFS. Finally
consider the IFS $\mathcal{F}= \{ \mathbb{P}^1; \, f_1,\, f_2\}$, where $\mathbb{P}^1$ is $1$-dimensional projective space, homeomorphic to the
circle, and which we denote by ${\mathbb{R}} \cup \{\infty\}$ (the reals
with an additional point at inf\/inity). The functions are
\begin{gather*}
f_1(x) = \frac{9x}{20-2x}, \qquad  f_2(x) = \frac{11x+9}{2x+18}.
\end{gather*}
The attractor $A$ of $\mathcal{F}$ is the interval $[0,1]$; the attractor $%
A^*$ of the dual IFS is ${\mathbb{R}}{\setminus} \{-9/2, 11/2\}$. Therefore
the basin of $A$ is $B = {\mathbb{R}}{\setminus} A^*= \{-9/2, 11/2\}$. It is
straightforward to check that the fast basin is $\widehat B = {\mathbb{R}}$.
Therefore, in this example, $B \subsetneq \widehat B$.

From the above examples we see that all three possibilities $\widehat B = B$,
$\widehat B \subsetneq B$, $B \subsetneq \widehat B$ can occur. However,
if the attractor $A$ has nonempty interior (in the topology of $\mathbb{X}$), then $B \subseteq \widehat{B}$ (see Corollary~\ref{cor:basins}).

Def\/ine the \textit{shift map} $S\colon \mathbb{I}_{0}\rightarrow\mathbb{I}_{0}$ by
\begin{gather*}
S(\iota_{1}\iota_{2}\cdots)=\iota_{2}\iota_{3}\cdots
\end{gather*}
for all $\iota\in\mathbb{I}_{0}$. Note that $S$ is continuous and maps $\mathbb{I}$ onto $\mathbb{I}$, $\mathbb{I}_+$ onto $\mathbb{I}_+$, and $\mathbb{I}_-$ onto~$\mathbb{I}_-$.

\begin{Definition}
A word $\theta \in \mathbb{I}_+$ is \textit{disjunctive} if every f\/inite
word is a subword of $\theta$. In fact, if~$\theta$ is disjunctive, then
every f\/inite word (in the alphabet~$I_+$) appears as a subword in~$\iota$
inf\/initely many times. Moreover, an equivalent def\/inition of a disjunctive
word is that its orbit under the shift map is dense in the code space.
Denote the interior of set~$S$ by~$S^{\circ}$ and the closure by~$\overline S$.
\end{Definition}

\begin{Theorem}\label{reversethm} Let $A$ be an attractor of an IFS $\mathcal{F}$, and let $\theta\in\mathbb{I}_{+}$ be disjunctive. Then $\pi(\theta)\in A^{\circ}$ if and only if $A^{\circ}\neq\varnothing$.
\end{Theorem}

\begin{proof}
Since $\theta\in\mathbb{I}_{+}$ is disjunctive, the orbit, $\left\{
S^{k}(\theta)\right\} _{k=0}^{\infty}$, of $\theta$ under the shift map $S\colon \mathbb{I}_{+}\mathbb{\rightarrow I}_{+}$ is dense in~$\mathbb{I}_{+}$.
Since $\pi$ is continuous, $\left\{ \pi_{\mathcal{F}}(S^{k}
(\theta))\right\} _{k=0}^{\infty}$ is dense in $A$. By way of contradiction
assume that $\pi_{\mathcal{F}}(\theta)\in\partial A$, the boundary of $A$.
Since each $f\in \mathcal{F}$ is a homeomorphism, $f(A^{\circ}) \subset
A^{\circ}$ or, equivalently, $\left ( f|_A\right )^{-1}(\partial A) \subset
\partial A$. Therefore $\pi(S^{k}(\theta))\in\partial A$ for all $k\in
\mathbb{N}_{0}$, which would mean that $A=\overline{\left\{
\pi(S^{k}(\theta))\right\} _{k=0}^{\infty}}\subset\overline{\partial A}%
=\partial A$ which implies $A^{\circ}=\varnothing$, which contradicts $%
A^{\circ} \neq \varnothing$.
\end{proof}

Theorem~\ref{interiorthm} below provides suf\/f\/icient conditions for a fractal
continuation of an attractor to contain the basin of the attractor. It
involves the equivalent notions of a full and a reversible word in $\mathbb{I}_+$. These concepts were introduced in~\cite{tiling}; the equivalence of
full and reversible is part of Theorem~\ref{thm:reversible}.

\begin{Definition}
If $A$ is an attractor of an IFS $\mathcal{F}$, call $\theta \in {\mathbb{I}}_+$ \textit{full} if there exists a nonempty compact set $A^{\prime}\subset
A^{\circ}$ such that, for any positive integer~$M$, there exist $n>m\geq M$ such
that
\begin{gather*}
f_{\theta_n}\circ f_{\theta_{n-1}} \circ \cdots \circ f_{\theta_{m+1}} (A)
\subset A^{\prime }.
\end{gather*}

Call $\theta \in \mathbb{I}_+$ \textit{reversible} w.r.t.\ an attractor $A$
and IFS $\mathcal{F}$ if there exists an $\omega=\omega_{1}\omega_{2}\cdots
\in\mathbb{I}_+$ such that $\omega$ is the address of some point in $A^{\circ}$
and, for every pair of positive integers $M$ and $L$, there is an integer $m\geq M$ such that
\begin{gather*}
\omega_{1}\omega_{2}\cdots \omega_{L} =
\theta_{m+L}\theta_{m+L-1}\cdots\theta_{m+1}.  \label{reveqn}
\end{gather*}
\end{Definition}

In some cases it is easy to check if a word is reversible. For example, if $\theta = \theta_1 \theta_2 \cdots \theta_k \theta_1 \theta_2 \cdots$ $\theta_k
\theta_1 \cdots$ is periodic and $\pi ( \theta_k \theta_{k-1} \cdots
\theta_1 \theta_k \theta_{k-1} \cdots \theta_1 \theta_{k} \cdots)$ lies in
the interior of $A$, then $\theta$ is reversible.

\begin{Theorem}[\protect{\cite[Theorem~3.7]{tiling}}] \label{thm:reversible} For an IFS $\mathcal{F}$, let $A$ be an attractor. With respect to~$A$ and~$\mathcal{F}$:

\begin{enumerate}\itemsep=0pt
\item[$1.$] There are inf\/initely many disjunctive words in $\mathbb{I}_+$ for $N\geq2$.

\item[$2.$] If $A^{\circ} \neq \varnothing$, then every disjunctive word is
reversible.

\item[$3.$] A word is reversible if and only if it is full.
\end{enumerate}
\end{Theorem}

\begin{Theorem}\label{interiorthm} Let $A$ be an attractor of an IFS $\mathcal{F}$ with
basin $B$. If $\theta\in\mathbb{I}_{+}$ is full/reversible, then
\begin{gather*}
B \subset {B}_{\theta}.
\end{gather*}
\end{Theorem}

\begin{proof}
Let $x\in B$. It suf\/f\/ices to show that $x \in f_{\theta_1}^{-1} \circ
f_{\theta_2}^{-1} \circ f_{\theta_3}^{-1} \circ \cdots \circ
f_{\theta_n}^{-1} (A)$ for some $n $, or equivalently
\begin{gather*}
f_{\theta_n} \circ f_{\theta_{n-1}} \circ \cdots \circ f_{\theta_1} (x) \in A.
\end{gather*}
By the def\/inition of attractor, for any $\epsilon > 0$ there is an $M_{\epsilon }$ such that if $m \geq M_{\epsilon }$, then
\begin{gather*}
F^m(x) \subset A_{\epsilon },
\end{gather*}
where $A_{\epsilon}$ is the open $\epsilon $-neighborhood of $A$. Because $\theta$ is assumed to be full, there exists a~compact set $A^{\prime }$ with
$A^{\prime}\subset A^{\circ}$ with the property that, for any $M$ there exist $n>m\geq M$ such that
\begin{gather*}
f_{\theta_n}\circ \cdots \circ f_{\theta_{m+1}}(A^{\prime}) \subset A^{\circ}.
\end{gather*}
This implies that there exists an $\epsilon _0$-neighborhood $A_{\epsilon
_0} $ of $A$ such that
\begin{gather*}
f_{\theta_n}\circ \cdots \circ f_{\theta_{m+1}}(A_{\epsilon _0}) \subset
A^{\circ},
\end{gather*}
for some $\epsilon _0>0$. If $M \geq M_{\epsilon _0}$, then
\begin{gather*}
f_{\theta_n}\circ \cdots \circ f_{\theta_{m+1}}\circ f_{\theta_m} \circ
\cdots \circ f_{\theta_1}(x) \in f_{\theta_n}\circ \cdots \circ
f_{\theta_{m+1}}(A_{\epsilon _0}) \subset A,
\end{gather*}
as required.
\end{proof}

\begin{Corollary}\label{cor:basins} Let $A$ be an attractor of an IFS $\mathcal{F}$ on a
complete metric space~$\mathbb{X}$, with basin~$B$ and fast basin $\widehat
B $. Then $B \subseteq \widehat B$ if and only if $A^{\circ} \neq \varnothing$. Moreover, if $\mathcal{F}$ is contractive and $A^{\circ} \neq \varnothing$,
then $B = \widehat B = \mathbb{X}$.
\end{Corollary}

\begin{proof}
Assume that $A^{\circ} \neq \varnothing$. Let $\theta$ be a disjunctive word;
by Theorem~\ref{thm:reversible} there are inf\/initely many. By the same
theorem $\theta$ is full and reversible. By Theorem~\ref{interiorthm} and
Proposition~\ref{prop:FB} we have $B \subseteq B_{\theta}\subseteq \widehat
B $. If $\mathcal{F}$ is contractive, then, as mentioned after the
def\/inition of contractive in Section~\ref{sec:IFS} (see also~\cite{Hutchinson}), it is well known that $B = \mathbb{X}$.

Conversely, assume that $A^{\circ} = \varnothing$. If $A$ is an attractor, the
fast basin of $A$ is
\begin{gather*}
\widehat B = \big\{ ( f_{\iota})^{-1}(A) \colon \iota \in {\mathbb{I}}_+^*\big\},
\end{gather*}
where ${\mathbb{I}}_+^*$ is the set of all f\/inite sequences in the alphabet~$I_+$. Since $(f_{\iota})^{-1}(A)$ has empty interior, $\widehat B$ is
nowhere dense by the Baire category theorem. Since~$B$ is open, clearly $B
\nsubseteq \widehat B$.
\end{proof}

\section{Symbolic IFSs and the extended coding map}\label{symbtheorysec}

The goal of this section is to extend the classic coding map given in
Def\/inition~\ref{def:CM}. We begin with two examples of ``symbolic IFSs''. It may be helpful for the reader, at this juncture, to recall the notation introduced around equation~\eqref{eq:notation} for $\mathbb{I}_0$, $\mathbb{I}$, $\mathbb{I}_+$, $\mathbb{I}_-$, $I$, $I_+$ and~$I_-$.

For
$n\in I$, def\/ine the inverse shift maps $\sigma_{n}\colon \mathbb{I}_0\rightarrow \mathbb{I}_0$, by
\begin{equation}  \label{eq:inShift0}
\sigma_{n}(\iota_{1}\iota_{2} \cdots)= n\iota_{1}\iota_{2} \cdots,
\end{equation}
and def\/ine $\sigma_{n}\colon \mathbb{I}\rightarrow\mathbb{I}$ by
\begin{gather*}  
\sigma_{n}(\iota_{1}\iota_{2} \cdots)=
\begin{cases}
n\iota_{1}\iota_{2} \cdots & \text{if} \ \ \iota_{1}\neq-n,\\
\iota_2\iota_3 \cdots &  \text{if} \ \ \iota_{1}=-n.
\end{cases}
\end{gather*}

\begin{Example}[first symbolic IFS]
\label{ex:Sym1} If
\begin{gather*}
Z = \{ \mathbb{I}_+; \, \sigma_n, \, n\in I_+ \},
\end{gather*}
then $Z$ is an IFS each of whose functions is a contraction with scaling
factor~$1/2$. The unique attractor is~$\mathbb{I}_+$. Not that this example
does not strictly f\/it the def\/inition of an IFS as given in Section~\ref{sec:IFS} because $\sigma_n$ is not a~homeomorphism as it is not surjective.
\end{Example}

\begin{Example}[second symbolic IFS]\label{ex:Sym2} If
\begin{gather*}
\widehat Z = \{\mathbb{I}; \,\sigma_{n},\, n\in I_{+} \},
\end{gather*}
then $\widehat Z$ is an IFS, for which it is easily verif\/ied that each
function $\sigma_n \colon  \mathbb{I}\rightarrow\mathbb{I}$ is bi-Lipshitz
with
\begin{gather*}
\frac12   d_{\mathbb{I}}(\iota,\omega)\leq d_{\mathbb{I}}\left( \sigma_{n}
(\iota),\sigma_{n}(\omega)\right) \leq 2  d_{\mathbb{I}}(\iota,\omega)
\end{gather*}
for all $\iota,\omega\in\mathbb{I}$. In particular, each $\sigma_n$ is a~homeomorphism with inverse $\sigma_{-n}\colon \mathbb{I\rightarrow}\mathbb{I}$.
The dual IFS is
\begin{gather*}
\widehat Z^* = \{ \mathbb{I}; \, \sigma_n, \, n \in I_- \}.
\end{gather*}
Although neither $\widehat Z$ nor $\widehat Z^*$ is contractive, both have
attractors as proved in Theorem~\ref{pairthm} below.
\end{Example}

In a metric space, the notation $\mathcal{B}(x,r)$ denotes the open ball of
radius~$r$ centered at~$x$.

\begin{Theorem}\label{pairthm}\quad
\begin{enumerate}\itemsep=0pt
\item[$1.$] The IFS $\widehat Z$ has a point-fibred attractor $\mathbb{I}_{+}$
with basin $\mathbb{I}\backslash \mathbb{I}_{-}$ and fast basin
\begin{gather*}  
\widehat {\mathbb{I}} := \left \{ \iota \in \mathbb{I }\colon \text{there
exists}\  k \in \mathbb{N} \  \text{such that} \ \iota_j \in
\begin{cases}
I_{-}   & \text{if} \ \ j < k \\
I_{+}   & \text{if} \ \  j\geq k
\end{cases}
\right \}.
\end{gather*}

\item[$2.$] The dual IFS $\widehat Z^*$ has point-fibred attractor $\mathbb{I}_{-}$
with basin $\mathbb{I}\backslash\mathbb{I}_{+}$ and fast basin
\begin{gather*}
{\widehat {\mathbb{I}}}^* := \left \{ \iota \in \mathbb{I }\colon \text{there exists} \  k \in \mathbb{N} \  \text{such that} \ \iota_j \in
\begin{cases}
I_{+}  & \text{if} \ \ j < k \\
I_{-}   & \text{if} \ \ j \geq k
\end{cases}
\right \}.
\end{gather*}
\end{enumerate}
\end{Theorem}

\begin{proof}
We will show that the IFS $\widehat Z$ has point-f\/ibred attractor $\mathbb{I}_{+}$ with basin $\mathbb{I} \backslash\mathbb{I}_{-}$ and dual repeller $\mathbb{I}_{-}$. The rest of the proof follows immediately, and is omitted.
The compact nonempty set $\mathbb{I}_{+}$ is contained in $\mathbb{I}$, and $\mathbb{I}_{+}=\widehat Z(\mathbb{I}_{+}):=\cup_{n\in I_{+}}\sigma_{n}(\mathbb{I}_{+})$. An open neighborhood of $\mathbb{I}_{+}$ is $\mathbb{I}\backslash\mathbb{I}_{-}$. If $\iota\in\mathbb{I}\backslash \mathbb{I}_{-}$,
then there is $K=K(\iota)\in\mathbb{N}$ such that $\iota_{K}\in I_{+}$, from
which it follows that $\alpha_{1}\in I_{+}$ for all $\alpha\in$ $\widehat
Z^{K-1}(\{\iota\})$. In turn, this implies
\begin{gather*}
d_{H}\big(\widehat Z^{K+j}(\{\iota\}),\mathbb{I}_{+}\big)\leq2^{-j-1},
\end{gather*}
for all $j\in\mathbb{N}$. This proves that $\big\{ \widehat
Z^{k}(\{\iota\})\big\} _{k\in\mathbb{N}}$ converges to $\mathbb{I}_{+}$
for all $\iota\in$ $\mathbb{I}\backslash\mathbb{I}_{-}$.

Now suppose that $C\in H(\mathbb{I})$ and $C\subset\mathbb{I}\backslash
\mathbb{I}_{-}$. Let $\iota\in C$, and let $K=K(\iota)$ be as above. Then
for all $\omega\in\mathcal{B}(\iota,2^{-K-1})$ we have $\omega_{K}\in I_{+}$
and hence
\begin{gather*}
d_{H}\big(\widehat Z^{K(\iota)+j}(\{\omega\}),\mathbb{I}_{+}\big)\leq2^{-j-1},
\end{gather*}
which shows that $\big\{ \widehat Z^{k}(\overline{\mathcal{B}
(\iota,2^{-K(\iota)-1})})\big\} _{k\in\mathbb{N}}$ converges to $\mathbb{I}
_{+}$. Using its compactness, $C$ can be covered by a f\/inite set of such
balls, and we conclude that $\big\{ \widehat Z^{k}(C)\big\} _{k\in\mathbb{N}}$ converges to~$\mathbb{I}_{+}$. On the other hand, $\mathbb{I}
_{-}\subset \widehat Z(\mathbb{I}_{-})$ so $\mathbb{I}_{-}\subset\widehat
Z^{k}(\mathbb{I}_{-})$ for all $k\in\mathbb{N}$, so $\big\{ \widehat Z^{k}(\mathbb{I}_{-})\big\} _{k\in\mathbb{N}}$ does not converge to~$\mathbb{I}_{+}$, so we conclude that the basin of the attractor $\mathbb{I}_{+}$
(w.r.t.~$\widehat Z$) is $\mathbb{I}\backslash\mathbb{I}_{-}$. Hence, by
def\/inition, the dual repeller is~$\mathbb{I}_{-}$.
\end{proof}

Theorem~\ref{HKtheorem} below is a generalization, along the lines of
Kieninger \cite[Section 4.2]{kieninger} and Hata \cite[Theorem 3.2]{hata},
to point-f\/ibred attractors, of ideas in Hutchinson's work \cite[2.1.(8),
p.~716; 3.1.(3), p.~724]{Hutchinson}. It asserts, in particular, a
semiconjugacy between~$\mathcal{F}$ acting on $A$ and the contractive IFS $Z
= \{\mathbb{I}_{+} ;\sigma_{n},n\in I_{+}\}$ of Example~\ref{ex:Sym1} acting
on its attractor,~$\mathbb{I}_{+}$.

\begin{Theorem}\label{HKtheorem} Let $A$ be a point-fibred attractor of the IFS $\mathcal{F}
=\{\mathbb{X};\,f_{n},\,n\in I_{+}\}$. The map $\pi\colon \mathbb{I}_{+}\rightarrow A$
from the attractor of~$Z$ to the attractor of~$\mathcal{F}$ given in
equation~\eqref{limiteqn} is uniformly continuous, and the diagram
\begin{gather*}
\xymatrix{
  \mathbb{I}_{+} \ar[r]^{\sigma_{n}}    \ar[d]_{\pi} & \mathbb{I}_{+}  \ar[d]^{\pi} \\
        A  \ar[r]_{f_{n}}  &   A   }
\end{gather*}
commutes for all $n\in I_{+}$.
\end{Theorem}

Theorem~\ref{HKtheorem} describes the behaviour of the functions of the IFS
on points of the attractor, not outside the attractor. The results below
describe the behaviour of the functions of the IFS on the fast basin of the
attractor and, in the case that the IFS is contractive with attractor with
nonempty interior, on all of the space $\mathbb{X}$. These new results,
Theorem~\ref{extenthm} and its corollaries, together with Theorem~\ref{HKtheorem}, provide a more complete symbolic description than the one
provided by Theorem \ref{HKtheorem} alone, of the structure and dynamics of
point-f\/ibred attractors of IFSs.

Extend the notation of equation~\eqref{eq:notation} as follows
\begin{gather*}
\mathbb{J}_{+}    =  \big\{\iota\in \mathbb{I} \colon \exists\,  K\in\mathbb{N}_{0}
\ \text{such that} \  S^{K}(\iota) \in\mathbb{I}_{+}\big\}, \\
\mathbb{J}_{-}    = \big\{\iota\in \mathbb{I}\colon \exists \, K\in\mathbb{N}_{0}
\ \text{such that} \  S^{K}(\iota) \in\mathbb{I}_{-}\big\},
\end{gather*}
and note that
\begin{gather*}
  \mathbb{I}_{+}  \subset \widehat {\mathbb I}
\subset\mathbb{J}_{+}\subset\overline{\mathbb{J}_{+} }=\mathbb{I}\subset
\mathbb{I}_{0},\qquad \text{and} \qquad \mathbb{I}_{-}   \subset {\widehat{\mathbb
I}}^* \subset\mathbb{J}_{-}
\subset\overline{\mathbb{J}_{-}}=\mathbb{I}\subset \mathbb{I}_{0}.
\end{gather*}
In particular, the spaces $\mathbb{J}_{\pm}$ are dense in $\mathbb{I}$.
Moreover, $\sigma_{n}(\mathbb{J}_{+})=\mathbb{J}_{+}$ and $\sigma_{n}(%
\mathbb{J}_{-})=\mathbb{J}_{-}$ for all $n\in I$. In addition, $\sigma_{n}|_{\mathbb{J}_{+}}\colon \mathbb{J}_{+}\rightarrow\mathbb{J}_{+}$ and $\sigma_{n}|_{\mathbb{J}_{-}}\colon \mathbb{J}_{-}\rightarrow\mathbb{J}_{-}$ are
homeomorphisms.

\begin{Definition}\label{welleqn-def}
Let $A$ be a point-f\/ibred attractor of an IFS $\mathcal{F}$. The \textit{extended coding map} $\pi \colon \mathbb{J}_{+}\rightarrow \mathbb{X}$ is
given by
\begin{equation}  \label{welleqn}
\pi(\iota)=\lim_{k\rightarrow\infty}f_{\iota|k}(b),
\end{equation}
for all $\iota\in\mathbb{J}_{+}$, where $b\in B(A)$. The limit in~\eqref{welleqn} exists and is independent of $b$ because, by the def\/inition
of $\mathbb{J}_{+}$, there is $K\in\mathbb{N}_{0}$ such that $S^{K}(\iota)\in
\mathbb{I}_{+}$, which means that we can write
\begin{gather*}
\pi(\iota) =f_{\iota|K}\big(\pi\big(S^{K}(\iota)\big)\big).
\end{gather*}
Clearly, the extended coding map agrees with the standard coding map
(Def\/inition~\ref{def:CM}) on~${\mathbb{I}}_+$. \textit{We will use the same
notation for both maps.} When the dual repeller~$A^*$ is a point-f\/ibred
attractor of the dual IFS ${\mathcal{F}}^*$, we denote the associated
extended coding map by ${\pi}^{\ast} \colon  {\mathbb{J}}_- \rightarrow
\mathbb{X}$.
\end{Definition}

Although $\mathbb{I}$ does not serve as a code space~-- in particular, a
result analogous to Theorem~\ref{HKtheorem} does not hold~-- the sets ${\mathbb{J}}_{\pm}$, that are dense in $\mathbb{I}$, do.

\begin{Theorem}\label{extenthm} Let $A$ be a point-fibred attractor of an IFS $\mathcal{F}$. The extended coding map \mbox{$\pi \colon \mathbb{J}_{+}\rightarrow \mathbb{X}$}
defined by~\eqref{welleqn} is continuous and agrees with the standard coding
map on~$\mathbb{I}_{+}$. The left diagram below is commutative for all $n
\in I = I_+\cup I_{-}$.

If $A^{\ast}$, the dual repeller of $A$, is a point-fibred attractor of the
dual IFS $\mathcal{F}^{\ast}$, then the right diagram below is commutative
for all $n\in I$:
\begin{gather*}
\xymatrix{
  \mathbb{J}_{+} \ar[r]^{\sigma_{n}}    \ar[d]_{\pi} & \mathbb{J}_{+}  \ar[d]^{\pi} \\
        \mathbb{X}  \ar[r]_{f_{n}}  &   \mathbb{X}   } \qquad
\xymatrix{
 \mathbb{J}_{-} \ar[r]^{\sigma_{n}}    \ar[d]_{\pi^\ast} & \mathbb{J}_{-}  \ar[d]^{\pi^\ast} \\
        \mathbb{X}  \ar[r]_{f_{n}}  &   \mathbb{X}   }
\end{gather*}
\end{Theorem}

\begin{proof}
That the extended coding map $\pi \colon \mathbb{J}_{+}\rightarrow \mathbb{X}$
agrees with the standard coding map on~$\mathbb{I}_{+}$ follows immediately
from Def\/initions~\ref{def:CM} and~\ref{welleqn-def}. That $\pi$ is continuous
follows from the def\/inition in Section~\ref{sec:IFS} of the metric on ${%
\mathbb{I}}_0$. Concerning the commuting diagram:
\begin{gather*}
\pi (\sigma_n(\iota)) = \lim_{k\rightarrow \infty} f_n\circ f_{\iota|k} (b)
)= f_n \big( \lim_{k\rightarrow \infty} f_{\iota|k} (b)\big) = f_n (\pi(\iota)).
\end{gather*}
The proof is similar for the dual.
\end{proof}

\begin{Example}
If $\mathcal{F}$ is as in Example~\ref{ex:schottky1}, then the ranges of
both~$\pi$ and~$\pi^{\ast} $ are contained in~$\Lambda$.
\end{Example}

\begin{Corollary}\label{theorem1} If $\mathcal{F}$ is an IFS with point-fibred attractor $A$
and fast basin $\widehat{B}$, then $\pi \colon \widehat{\mathbb{I}}\rightarrow
\widehat B$ is surjective and the left diagram below is commutative for all $n\in I_{-}$.

Let ${\widehat B}^*$ denote the fast basin of the dual IFS $\mathcal{F}^*$.
If $A^{\ast}=\mathbb{X}\backslash B(A)$ is a point-fibred attractor of the
dual IFS $\mathcal{F}^{\ast}$, then ${\pi}^* \colon {\widehat{\mathbb{I}}}
^*\rightarrow {\widehat B}^*$ is surjective and the right diagram below is
commutative for all $n\in I_{+}$.
\begin{gather*}
\xymatrix{
  \widehat{\mathbb{I}} \ar[r]^{\sigma_{n}}    \ar[d]_{\pi} & \widehat{\mathbb{I}}  \ar[d]^{\pi} \\
        \widehat{B}  \ar[r]_{f_{n}}  &   \widehat{B}   } \qquad
\xymatrix{
 \widehat{\mathbb{I}}^{\ast} \ar[r]^{\sigma_{n}}    \ar[d]_{\pi^\ast} & \widehat{\mathbb{I}}^{\ast}  \ar[d]^{\pi^\ast} \\
        {\widehat B}^{\ast}  \ar[r]_{f_{n}}  &   {\widehat B}^{\ast}   }
\end{gather*}
\end{Corollary}

\begin{proof}
This is a corollary of Theorem \ref{extenthm}. We simply note that $%
\sigma_{-n}(\widehat{\mathbb{I}})\subset\widehat{\mathbb{I}}$ and $%
\sigma_{n}(\widehat{\mathbb{I}}^{\ast})\subset\widehat{\mathbb{I}}^{\ast}$
for all $n\in I_{+}$.
\end{proof}

The following corollary provides a model for attractor-repeller pairs, for
example the loxodromic M\"{o}bius case discussed in \cite{Vince}.

\begin{Corollary}
\label{dualcasethm} Let $\mathcal{F}$ be an IFS with point-fibred attractor $%
A$ on a complete metric space $X$ such that $A^{\circ}\neq\varnothing$. If $%
\mathcal{F}$ is contractive, then ${\pi}(\widehat{\mathbb{I}} )=\mathbb{X}$
and the left diagram below commutes for all $n\in I_{-}$. If $\mathbb{X}$ is
compact, let the dual repeller $A^{\ast}$ be a point-fibred attractor of $%
\mathcal{F}^{\ast}$ with $(A^{\ast})^{\circ }\neq\varnothing$. Then ${\pi}(%
\widehat{\mathbb{I}} )={\pi}^*(\widehat{\mathbb{I}}^*)=\mathbb{X}$; the left
diagram commutes for all $n\in I_{-}$, and the right diagram commutes for
all $n\in I_{+}$.
\begin{gather*}
\xymatrix{
  \widehat{\mathbb{I}} \ar[r]^{\sigma_{n}}    \ar[d]_{\pi} & \widehat{\mathbb{I}}  \ar[d]^{\pi} \\
        \mathbb{X}  \ar[r]_{f_{n}}  &   \mathbb{X}   } \qquad
\xymatrix{
 \widehat{\mathbb{I}}^{\ast} \ar[r]^{\sigma_{n}}    \ar[d]_{\pi^\ast} & \widehat{\mathbb{I}}^{\ast}  \ar[d]^{\pi^\ast} \\
        \mathbb{X}  \ar[r]_{f_{n}}  &   \mathbb{X}   }
\end{gather*}
\end{Corollary}

\begin{proof}
The result follows from Corollary~\ref{theorem1} once it is shown that the
fast basin of $A$ is ${\mathbb{X}}$ and the fast basin of $A^*$ is also ${%
\mathbb{X}}$. In the case of a contractive IFS with $A^{\circ} \neq \varnothing
$, the fast basin and the basin coincide and both equal $\mathbb{X}$; see
Section~\ref{sec:IFS}. The dual repeller in this case is empty.

In the $\mathbb{X}$ compact case, we will show that the fast basin of $A^*$
is ${\mathbb{X}}$; a similar argument shows that the fast basin of $A$ is ${\mathbb{X}}$. It is known \cite[Theorem~5.2]{BVaustMS} that the basin $B(A^*) = {\mathbb{X}} {\setminus} A$. So by Proposition~\ref{prop:FB} and
Theorems~\ref{interiorthm} and \ref{thm:reversible}, if $x \notin A$, then $x $ is in the fast basin of~$A^*$. Lastly, if $x \in A$, let $y\neq x$ be
any other point of~$A$. If~$\iota$ is any address of $y$, then $\pi(\iota) =
\lim\limits_{k\rightarrow \infty} f_{\iota|k}(A)$. Therefore there is a $k$ such
that $x \notin f_{\iota|k}(A)$, which implies that
\begin{gather*}
z := f_{\iota_k}^{-1}\circ f_{\iota_{k-1}}^{-1} \circ \cdots \circ
f_{\iota_1}^{-1} (y) \notin A.
\end{gather*}
Since $z$ is in the basin of $A^*$ and $(A^*)^\circ \neq \varnothing$, there is $\omega_1 \omega_2 \cdots \omega_j$ for some $j$ such that $%
f_{\omega_1}^{-1}\circ f_{\omega_{2}}^{-1} \circ \cdots \circ
f_{\omega_j}^{-1} (z) \in A$. Hence $x$ lies in the fast basin of~$A^*$.
\end{proof}

Theorem \ref{extenthm} and its corollaries have many consequences including
the following: (i)~the feasibility of continuous assignment of addresses to
points in~$\widehat{B}(A)$; (ii)~symbolic dynamics on~$\widehat{\mathbb{I}}$
are semiconjugate to corresponding dynamics on $\widehat{B}(A)$; (iii)~if $A$
has nonempty interior as a~subset of~$\mathbb{X}$, then we obtain addresses
and dynamics on the basin~$B(A)$; (iv)~description of the relationship
between attractors and dual repellers; (v)~extension of fractal
homeomorphisms, generalizing results in~\cite{monthly} and~\cite{tiling};
(vi)~address structures for continuations, branched fractal manifolds and
tilings, thereby clarifying and extending ideas in~\cite{tiling}.

\section{Fractal manifold}\label{leafsec}

Throughout this section, $A$ is a point-f\/ibred attractor of an IFS $\mathcal{F}$ on a complete metric space~$\mathbb{X}$. The branched fractal manifold
generated by~$(\mathcal{F}, A)$ is constructed in this section. Properties
of the branched fractal manifold are the subject of Theorem~\ref{metricthm}.

\begin{Definition}
For $\iota\in\widehat{\mathbb{I}}$, let
\begin{gather*}
k_{\iota}=\min\big\{k\colon {\pi} \big(S^{k}(\iota)\big)\in A\big\}\in\mathbb{N}_{0}.
\end{gather*}
Def\/ine
\begin{gather*}
\lbrack\iota]=\iota|k_{\iota}\qquad \text{and} \qquad \left\langle \iota\right\rangle
=S^{k_{\iota}}(\iota),
\end{gather*}
to be the \textit{integer part} and the \textit{fractional part} of the
address $\iota = \iota_1 \iota_2 \cdots \iota_{k_{\iota}} \cdots$ (w.r.t.~$\mathcal{F}$ and~$A$) respectively.
\end{Definition}

\begin{Example}[integer and fractional parts]
Consider the case $\mathcal{F} =\{\mathbb{R};\,f_{1},\,f_{2}\}$, where $f_{1}(x)=\frac{1}{2}x$ and $f_{2} (x)=\frac{1}{2}x+\frac{1}{2}$. The unique attractor
is $A=[0,1]$. We will determine the integer and fractional parts of the two
elements of $\widehat{\mathbb{I}}$, namely $-1 \,-1\,-1\, \overline 2$ and $-1 \, -1\, -2\, 1 \, \overline 2$. A ``bar'' over a sequence means inf\/inite
repetition; in particular $\overline 2 = 222 \cdots$. Noting that $%
\pi(\overline 2) = 1 \in {\mathbb{R}}$ and $\pi(-1 \, \overline{2})= 2 \in {\mathbb{R}}$ we have
\begin{gather*}
  [ -1\, -1\, -1\, \overline{2} ]  = -1\, -1\, -1, \qquad
\langle -1\, -1\, -1\, \overline{2} \rangle  = \overline 2  .
\end{gather*}
Noting that $\pi( -1 \, -1 \, -2 \, 1\, \overline{2})= 0 \in {\mathbb{R}}$,
we have
\begin{gather*}
  \left [ -1 \, -1 \, -2 \, 1 \, \overline{2} \right ]  =
\varnothing, \qquad \langle -1\, -1 \, -2\, 1\, \overline{2} \rangle  = -1\, -1 \,
-2\, 1\, \overline{2}.
\end{gather*}
\end{Example}

Def\/ine an equivalence relation $\sim$ on $\widehat{\mathbb{I}}$ by $\iota
\sim \omega$ if and only if
\begin{enumerate}\itemsep=0pt
\item[1)] $[\iota] = [\omega]$, and

\item[2)] $\pi (\iota) = \pi (\omega)$,
\end{enumerate}
and let
\begin{gather*}
\mathbb{L }= \widehat{\mathbb{I}}/{\sim}.
\end{gather*}
For $\iota\in\widehat{\mathbb{I}}$, let $\widetilde{\iota}$ denote the
equivalence class to which $\iota$ belongs. Maps $\widetilde{\pi }\colon \widehat{%
\mathbb{I}}\rightarrow\mathbb{L}$ and $\widehat{\pi}\colon \mathbb{L\rightarrow}
\mathbb{X}$ are well def\/ined by $\widetilde{\pi} (\iota)=\widetilde{\iota}$
and $\widehat{\pi}(\widetilde{\iota})={\pi} (\iota)$, respectively. The maps
$\widehat \pi$ and $\widetilde{\pi}$ will be called \textit{projection maps}.

For $\alpha:=\widetilde{\iota}\in\mathbb{L}$, let
\begin{gather*}
k_{\alpha}:=k_{\iota}, \qquad   [\alpha]:=[\iota],\qquad
\left\langle \alpha \right\rangle :=\widetilde{\left\langle
\iota\right\rangle}.
\end{gather*}
These entities do not depend on the representative $\iota$ of the
equivalence class~$\alpha$. It is useful to think of~$\mathbb{L}$ in terms
of leaves, panicles, and sheets.

\begin{Definition}
If $\theta$ is a f\/inite word in the alphabet~$I_-$, then
\begin{gather*}
l(\theta) = \{ \alpha \in \mathbb{L}\colon  [\alpha ] = \theta\}
\end{gather*}
will be called a \textit{leaf} of $\mathbb{L}$, and
\begin{gather*}
p(\theta)=l(\varnothing) \cup l(\theta_{1})\cup l(\theta_{1}\theta_{2})\cup
\cdots \cup l(\theta)
\end{gather*}
will be called a \textit{panicle} of~$\mathbb{L}$. Because of the special
importance of the case $\theta=\varnothing$, def\/ine
\begin{gather*}
\mathbb{K} := l (\varnothing) = p(\varnothing).
\end{gather*}
If $\theta \in {\mathbb{I}}_{-}$, we will refer to the ``inf\/inite panicle''
\begin{gather*}
\mathbb{L}_{\theta} := \bigcup_{k=0}^{\infty} l(\theta|k)
\end{gather*}
as a \textit{sheet} of $\mathbb{L}$.
\end{Definition}

Def\/ine a metric on $\mathbb{L}$ as follows. That it is a metric is part of
Theorem~\ref{metricthm} below. For $\alpha, \beta \in \mathbb{L}$ let
\begin{gather*}
  K  := \max \left \{ k \colon  [\alpha]|k = [\beta]|k \right \},\qquad
  [\alpha,\beta]  := [\alpha]|K = [\beta]|K, \qquad  p(\alpha, \beta)  := p(
[\alpha, \beta]),
\end{gather*}
and
\begin{gather}
d_{\mathbb{L}}(\alpha,\beta)=\min  \big\{d_{\mathbb{X}}(\widehat{\pi}(\alpha ),\widehat{\pi}(\iota))+d_{\mathbb{X}}(\widehat{\pi}(\iota),\widehat{\pi }(\beta))\colon \iota\in p(\alpha,\beta)\big\}.  \label{sevenptone}
\end{gather}
Because $p(\alpha,\beta) = f_{ [\alpha,\beta] }(A)$, where the notation of
equation~\ref{eq:composition} is used, it is equivalent to def\/ine
\begin{gather}  \label{eq:metric}
d_{\mathbb{L}} (\alpha, \beta) = \min \big \{ d_{\mathbb{X}} (
\widehat{\pi}(\alpha), x  ) + d_{\mathbb{X}}  (x, \widehat{\pi}
(\beta)  ) \colon x \in f_{ [\alpha,\beta] }(A) \big \}.
\end{gather}

The metric space $({\mathbb{L}},d_{\mathbb{L}})$ will be called the \textit{branched fractal manifold} or just the \textit{f-manifold} generated by $(\mathcal{F},A)$. A description of the f-manifold is the intent of the
following theorem. Shift maps $\widetilde \sigma_n \colon \mathbb{L} \rightarrow \mathbb{L} $ are well-def\/ined by
\begin{gather*}
\widetilde \sigma_n (\widetilde \iota) = \widetilde{\sigma_n(\iota)}.
\end{gather*}
These maps occur in statement~(9).

\begin{Theorem}
\label{metricthm} Let $A$ be a point-fibred attractor of an IFS $\mathcal{F}$
with fast basin $\widehat{B}\subset\mathbb{X}$. Let $(\mathbb{L},d_{\mathbb{L}})$ be the f-manifold generated by $(\mathcal{F},A)$. Let $\widehat{\pi}\colon
\mathbb{L\rightarrow} \mathbb{X}$ and $\widetilde{\pi }\colon \widehat{\mathbb{I}}%
\mathbb{\rightarrow L}$ be the associated projection maps. The following
statements hold.
\begin{enumerate}\itemsep=0pt

\item[$1.$] The minimum in equation~\eqref{sevenptone} is achieved, and~$d_{\mathbb{L}}$ is a metric on~$\mathbb{L}$.

\item[$2.$] The set of all leaves form a partition of~$\mathbb{L}$, and each leaf
is nonempty.

\item[$3.$] The image of $\mathbb{L}$ under $\widehat{\pi}$ is the fast basin,
i.e., $\widehat \pi (\mathbb{L}) = \widehat B$.

\item[$4.$] For any $\theta\in \mathbb{I}_{-}$, the projection $\widehat \pi$ maps
the sheet $\mathbb{L}_{\theta}$ isometrically onto the continuation $B_{-\theta} \subseteq \widehat B \subseteq \mathbb{X}$:
\begin{gather*}
d_{\mathbb{X}}(\widehat{\pi}(\alpha),\widehat{\pi}(\beta)) = d_{\mathbb{L}}(\alpha,\beta).
\end{gather*}
for all $\alpha,\beta\in\mathbb{L}_{\theta}$.

\item[$5.$] Each panicle is homeomorphic to $A$.

\item[$6.$] The projection maps $\widetilde{\pi}\colon \widehat{\mathbb{I}}\rightarrow \mathbb{L}$ and $\widehat{\pi}\colon \mathbb{L}\rightarrow\widehat{B}$, are
continuous; in particular
\begin{gather*}
d_{\mathbb{X}}(\widehat{\pi}(\alpha),\widehat{\pi}(\beta))\leq d_{\mathbb{L}}(\alpha,\beta).
\end{gather*}
for all $\alpha,\beta\in\mathbb{L}$.

\item[$7.$] If $A$ is pathwise connected, then $\mathbb{L}$ is pathwise connected. %

\item[$8.$] If the metric spaces $(\mathbb{X},d_{\mathbb{X}})$ and $(\mathbb{X}
,d_{\mathbb{X}}^{\prime})$ have the same topology, then $\left( \mathbb{L},d_{\mathbb{L}}\right) $ and $(\mathbb{L},d_{\mathbb{L}}^{\prime})$ have the
same topology, where the metric $d_{\mathbb{L}}^{\prime}$ is defined by
equation~\eqref{sevenptone} with $d_{\mathbb{X}}$ replaced by $d_{\mathbb{X}}^{\prime}$.

\item[$9.$] With $\sigma_{n}\colon \widehat{\mathbb{I}}\rightarrow\widehat{\mathbb{I}}$
as defined in Section~{\rm \ref{symbtheorysec}} and $\widetilde{\sigma}_{n}\colon
\mathbb{L}\rightarrow \mathbb{L}$ as defined above, both $\sigma_n$ and $\widetilde{\sigma}_n$ are continuous and injective, for all $n\in I_{-}$, and we have
the commutative diagram,
\begin{gather*}
\xymatrix{
  \widehat{\mathbb{I}} \ar[r]^{\sigma_{n}}    \ar[d]_{\widetilde{\pi}} & \widehat{\mathbb{I}}  \ar[d]^{\widetilde{\pi}} \\
        \mathbb{L}  \ar[r]_{\widetilde{\sigma}_{n}}\ar[d]_{\widehat{\pi}}  &   \mathbb{L}  \ar[d]^{\widehat{\pi}} \\
  \widehat{B}  \ar[r]_{f_{n}}  &   \widehat{B}   }
\end{gather*}
\end{enumerate}
\end{Theorem}

\begin{proof}
The proofs of (2), (3), (9) are straightforward.

(1) The set $f_{[\alpha,\beta]}(A)$ in \eqref{eq:metric} is compact because $%
A$ is compact and $f_{[\alpha,\beta]}$ is continuous. For f\/ixed $\alpha
,\beta$, the function $d_{\mathbb{X}}(\widehat{\pi}(\alpha),x)+d_{\mathbb{X}
}(x,\widehat{\pi}(\beta))$ is continuous in $x\in f_{[\alpha,\beta]}(A)$.
Hence the minimum is achieved at some $x^{\ast}\in f_{[\alpha,\beta]}(A)$.

Concerning the metric, we establish only the triangle inequality:
\begin{gather*}
d_{\mathbb{L}}(\alpha ,\gamma)+d_{\mathbb{L}}(\gamma,\beta)=d_{\mathbb{X}}(\widehat{\pi} (\alpha),x_{1})+d_{\mathbb{X}}(x_{1},\widehat{\pi}(\gamma))+d_{\mathbb{X} }(\widehat{\pi}(\gamma),x_{2})+d_{\mathbb{X}}(x_{2},\widehat{\pi}
(\beta)),
\end{gather*}
where $x_{1}\in f_{[\alpha,\gamma]}(A)$ minimizes $d_{\mathbb{X}}(\widehat {\pi}(\alpha),x_{1})+d_{\mathbb{X}}(x_{1},\widehat{\pi}
(\gamma))$ and $x_{2}\in f_{[\beta,\gamma]}(A)$ minimizes \linebreak $d_{\mathbb{X}}(\widehat{\pi }(\gamma),x_{2})+d_{\mathbb{X}}(x_{2},\widehat{\pi}(\beta))$.
We have either $f_{[\alpha,\gamma]}(A)\subset f_{[\alpha,\beta]}(A)$ or $f_{[\beta,\gamma ]}(A)\subset f_{[\alpha,\beta]}(A)$. Without loss of
ge\-ne\-rality assume that $f_{[\beta,\gamma ]}(A)\subset f_{[\alpha,\beta]}(A)$. Then
\begin{gather*}
d_{\mathbb{L}}(\alpha,\gamma)+d_{\mathbb{L}}(\gamma,\beta)   =d_{\mathbb{X}
}(\widehat{\pi}(\alpha),x_{1})+d_{\mathbb{X}}(x_{1},\widehat{\pi}
(\gamma))+d_{\mathbb{X}}(\widehat{\pi}(\gamma),x_{2})+d_{\mathbb{X}} (x_{2},
\widehat{\pi}(\beta)) \\
\hphantom{d_{\mathbb{L}}(\alpha,\gamma)+d_{\mathbb{L}}(\gamma,\beta)}{}
 \geq d_{\mathbb{X}}(\widehat{\pi}(\alpha),x_{1})+d_{\mathbb{X}}(x_{1}
,x_{2})+d_{\mathbb{X}}(x_{2},\widehat{\pi}(\beta)) \\
\hphantom{d_{\mathbb{L}}(\alpha,\gamma)+d_{\mathbb{L}}(\gamma,\beta)}{}
 \geq d_{\mathbb{X}}(\widehat{\pi}(\alpha),x_{2})+d_{\mathbb{X}}(x_{2},
\widehat{\pi}(\beta))\text{ (where }x_{2}\in f_{[\beta,\gamma ]}(A)\text{)}
\\
\hphantom{d_{\mathbb{L}}(\alpha,\gamma)+d_{\mathbb{L}}(\gamma,\beta)}{}
 \geq\min\{d_{\mathbb{X}}(\widehat{\pi}(\alpha),x)+d_{\mathbb{X}}(x,
\widehat{\pi}(\beta))\colon x\in f_{[\alpha,\beta]}(A)\} \\
\hphantom{d_{\mathbb{L}}(\alpha,\gamma)+d_{\mathbb{L}}(\gamma,\beta)}{}
 =d_{\mathbb{L}}(\alpha,\beta).
\end{gather*}

(4) Assume that $\alpha, \beta \in \mathbb{L}$ lie in the same sheet, and
therefore lie in a common panicle, say $p([\beta])$. Then $[\alpha,\beta] =
[\alpha]$ and $f_{[\alpha,\beta]}(A)=f_{[\alpha]}(A)$, which implies that $
\widehat{\pi}(\alpha)\in f_{[\alpha]}(A)$. By the triangle inequality we
have $d_{\mathbb{L}}(\alpha,\beta)=\min\{d_{\mathbb{X}}(\widehat{\pi}
(\alpha),x)+d_{\mathbb{X}}(x,\widehat{\pi}(\beta))\colon x\in f_{[\alpha,\beta
]}(A)\}=d _{\mathbb{X}}(\widehat{\pi}(\alpha),\widehat{\pi}(\beta))$. That the image of $\mathbb{L}_{\theta}$ under $\widehat \pi$ is $B_{-\theta}$ follows from the def\/initions of $\widehat \pi$ and of the
continuation~$B_{-\theta}$.

(5) It follows from the def\/initions that $\widehat \pi$ maps the panicle $p(\theta)$ onto~$f_{\theta}(A)$. The result now follows from that facts that
$f_{\theta}$ is a homeomorphism and that, by part~(4), the mapping $\widehat
\pi$ restricted to~$p(\theta)$ is an isometry.

(6) The continuity of $\widehat{\pi}\colon \mathbb{L}\rightarrow\mathbb{X}$
follows from the inequality $d_{\mathbb{X}}(\widehat{\pi }(\alpha),\widehat{\pi}(\beta)) \leq d_{\mathbb{L}}(\alpha,\beta)$, which follows from the
triangle inequality.

Concerning the continuity of $\widetilde{\pi}\colon \widehat{\mathbb{I}}\rightarrow
\mathbb{L}$, let $\iota\in\widehat{\mathbb{I}}$ be given. If $\omega\in
\widehat{\mathbb{I}}$ and $\omega|K=\iota|K$ for $K$ suf\/f\/iciently large,
that is, if $d_{\mathbb{I}}(\iota,\omega)$ is suf\/f\/iciently small, then $\widetilde \omega$ and $\widetilde \iota$ lie the same leaf. By part~(4) $d_{\mathbb{L}}(\widetilde{\pi}(\iota),\widetilde{\pi}(\omega)) = d_{\mathbb{X}}
(\widehat{\pi} (\widetilde {\iota}),\widehat {\pi} (\widetilde {\omega}))=
d_{\mathbb{X}}(\widehat \pi(\iota),\widehat \pi(\omega))$. The continuity of
$\widetilde{\pi}\colon \widehat{\mathbb{I}}\rightarrow\mathbb{L}$ then follows
from the continuity of $\widehat \pi\colon \widehat {\mathbb{I}}\rightarrow\mathbb{%
X}$ by Theorem~\ref{extenthm}.

(7) Assume that $A$ is pathwise connected. Given any two points $\alpha,
\beta \in \mathbb{L}$, the leaf $\mathbb{K}$ is contained in the
intersection of the sheets containing~$\alpha$ and~$\beta$. It is therefore
suf\/f\/icient to prove that any sheet is pathwise connected; so assume that~$\alpha$ and~$\beta$ lie in the same sheet, which implies that they lie in
the same panicle $p(\theta)$ for some f\/inite word $\theta$. But by~(5), this
panicle is homeomorphic to $A$.

(8) Two metrics $d_{\mathbb{X}}$ and $d_{\mathbb{X}}^{\prime}$ on $\mathbb{X}
$ generate the same topology if, given any ball $\mathcal{B}_{(\mathbb{X}
,d_{\mathbb{X}})}(x,r)$, \mbox{$r>0$}, there is a positive radius $r^{\prime}$
such that $\mathcal{B}_{(\mathbb{X},d_{\mathbb{X}}^{\prime})}(x,r^{\prime})
\subset\mathcal{B}_{(\mathbb{X},d_{\mathbb{X}})}(x,r)$ and, given any \linebreak $
\mathcal{B}_{(\mathbb{X},d_{\mathbb{X}}^{\prime})}(x,r^{\prime})$, there is
a positive radius $r$ such that $\mathcal{B}_{(\mathbb{X},d_{\mathbb{X}}
)}(x,r)\subset\mathcal{B}_{(\mathbb{X},d_{\mathbb{X}}^{\prime})}(x,r^{\prime
})$. In the present situation, let the metrics on $\mathbb{L}$ def\/ined using
$d_{\mathbb{X}}$ and $d_{\mathbb{X}}^{\prime}$ be denoted by $d_{\mathbb{L}}$
and $d_{\mathbb{L}}^{\prime}$. Let $r>0$ be given and $\alpha\in \mathbb{L}$
be given. We will show that there is $r^{\prime}>0$ so that $\mathcal{B}_{(
\mathbb{L},d_{\mathbb{L}}^{\prime})}(\alpha,r^{\prime} )\subset\mathcal{B}_{(
\mathbb{L},d_{\mathbb{L}})}(\alpha,r)$. The rest of the argument is then
obtained by switching the roles of the key players.

If $\alpha$ and $\beta$ are contained in the same sheet, then by~(4), we
have $d_{\mathbb{L}}(\alpha,\beta )=d_{\mathbb{X}}(\widehat{\pi}(\alpha),\widehat{\pi}(\beta))$, so that
\begin{gather*}
\mathcal{B}_{(\mathbb{L},d_{\mathbb{L}})}(\alpha,r) =\big\{\beta\in\mathbb{L}\colon
\widehat{\pi}(\beta)\in \mathcal{B}_{(\mathbb{X},d_{\mathbb{X}})}(\widehat{\pi}(\alpha),r)\big\}.
\end{gather*}
Now choose $r^{\prime}>0$ so small that
\begin{gather*}
\mathcal{B}_{(\mathbb{X},d_{\mathbb{X}}^{\prime})}(\widehat{\pi}
(\alpha),r^{\prime })\subset\mathcal{B}_{(\mathbb{X},d_{\mathbb{X}})}(
\widehat{\pi}(\alpha),r).
\end{gather*}
Then, if $r^{\prime}$ is suf\/f\/iciently small,
\begin{gather*}
\mathcal{B}_{(\mathbb{L},d_{\mathbb{L}})}(\alpha,r') =
\{\beta\in\mathbb{L}\colon \widehat{\pi}(\beta)\in
\mathcal{B}_{(\mathbb{X},d_{\mathbb{X}})}(\widehat{\pi}(\alpha),r')\}\\
\hphantom{\mathcal{B}_{(\mathbb{L},d_{\mathbb{L}})}(\alpha,r')}{}
 \subset \{\beta\in\mathbb{L}\colon \widehat{\pi}(\beta)\in
\mathcal{B}_{(\mathbb{X},d_{\mathbb{X}})}(\widehat{\pi}(\alpha),r)\}
=\mathcal{B}_{(\mathbb{L},d_{\mathbb{L}})}(\alpha,r) .
\end{gather*}

It now suf\/f\/ices to show that if $\alpha$ is f\/ixed and $\beta$ is
suf\/f\/iciently close to $\alpha$, then $\alpha$ and $\beta$ lie on the same
sheet. Assume that $\{\beta_n\}$ is a sequence of points in~$\mathbb{L}$
that converges to $\alpha$, but $\beta_n$ is not in the same sheet as $\alpha
$ for every $n$. By the def\/inition of the distance on the f-manifold, this
implies that there is a sequence $\{\gamma_n\}$ such that $\gamma_n \in
p(\alpha,\beta_n)$ and $\widehat \pi (\gamma_n)$ converges to $\widehat \pi
(\alpha)$. Because $\alpha$ is f\/ixed, the set $\{p(\alpha,\beta_n)\}$ has at
most f\/initely many distinct elements. So there is an inf\/inite subsequence $\{\gamma^{\prime }_n\}$ of $\{\gamma_n\}$ and the corresponding subsequence $\{\beta^{\prime }_n\}$ of $\{\beta_n\}$ such that
\begin{gather*}
\gamma^{\prime }_n \in p(\alpha, \beta^{\prime }_n) = p(\alpha,\beta_0)
\end{gather*}
for some f\/ixed $\beta_0$. Since $\alpha, \gamma^{\prime }_n \in p(\alpha)$
for all $n$ and $p(\alpha)$ is compact by statement~(5), it follows from the
fact that $\widehat \pi (\gamma^{\prime }_n)$ converges to $\widehat \pi
(\alpha)$ and statement~(4) that $\{\gamma_n\}$ converges to~$\alpha$. But
from the facts that $p(\alpha, \beta_0)$ is compact and $\gamma_n^{\prime
}\in p(\alpha, \beta_0)$, it follows that $\alpha \in p(\alpha,\beta_0)$.
Therefore, $\alpha$ and $\beta^{\prime }_n$ are in the same panicle, and
hence the same sheet, for all~$n$, a contradiction.
\end{proof}

\begin{Theorem}
\label{sevenpttwothm} Let $A$ be a point-fibred attractor of an IFS $\mathcal{F}$. For the f-manifold of $(\mathcal{F},A)$ there are at most $N +1
$ topologically distinct leaves. Specifically, the leaf $l(\varnothing)=\mathbb{K}$ is homeomorphic to $A$, and the leaf $l(\theta)$ for $\theta\neq\varnothing$ is homeomorphic to $A\backslash f_{i_{\theta}}(A)$
where~$-i_{\theta}$ is the last component of~$\theta$.
\end{Theorem}

\begin{proof}
To simplify notation, when the intent is clear, we make no distiction
between $\iota\in\widehat{\mathbb{I}}$ and its equivalence class in $\mathbb{L}$. For $\alpha\in\mathbb{L}$, note that ${\pi}(\langle \alpha\rangle)\in A$. For a leaf $l(\theta)$, let $\langle l(\theta)\rangle =\{\langle \alpha
\rangle\colon \alpha\in l(\theta)\}$. By def\/inition, $\iota\in\langle
l(\theta)\rangle$ if and only if $\pi(\iota)\in A$ and $f_{S^{i}(\iota)}(%
\pi(\iota))\notin A$ for $i=0,1,2,\dots ,k_{\iota}-1$; the last condition here
is equivalent to $f_{i_{\theta}}(\pi(\iota))\notin A$ where $-i_{\theta}$ is
the last entry in $\theta$. Therefore ${\pi }(\langle l(\theta)
\rangle)=A\backslash f_{i_{\theta}}(A)$. This implies that $\widehat{\pi}(l(\theta))=f_{\theta}(\pi(\langle l(\theta)\rangle))=f_{\theta
}(A\backslash f_{i_{\theta}}(A))$. Since $\widehat \pi$ restricted to $l(\theta)$ is a homeomorphism by part (4) of Theorem~\ref{metricthm}, and $f_{\theta}$ is a homeomorphism, we have $l(\theta)$ is homeomorphic to $A\backslash f_{i_{\theta}}(A)$.
\end{proof}

\begin{Corollary}
\label{cor:leaf} Let $A$ be a point-fibred attractor of an IFS $\mathcal{F}$. The closures of any two distinct leaves of the f-manifold are
non-overlapping.
\end{Corollary}

\begin{proof}
This follows from part (5) of Theorem~\ref{metricthm}, Theorem~\ref{sevenpttwothm}, and the fact that the boundary of a compact set has empty
interior.
\end{proof}

\begin{Example}
The attractor of the IFS $\{ {\mathbb{R}} ; \, f_1, \,f_2 \}$ where $f_1(x) =
\frac13 x$ and $f_2(x) = \frac13 x + \frac23$, is the Cantor set~$\mathcal{C}
$. The f-manifold~$\mathbb{L}$ is the disjoint union of its leaves, each of
which is homeomorphic to $\mathcal{C}$. The image of each leaf under the
projection map $\widehat \pi$ is related by a similitude to~$\mathcal{C}$,
with scaling factor equal to~$3^k$ for some $k \in {\mathbb{N}}_0$.
\end{Example}

\begin{figure}[t]
 \centering
\includegraphics[width=12cm, keepaspectratio]{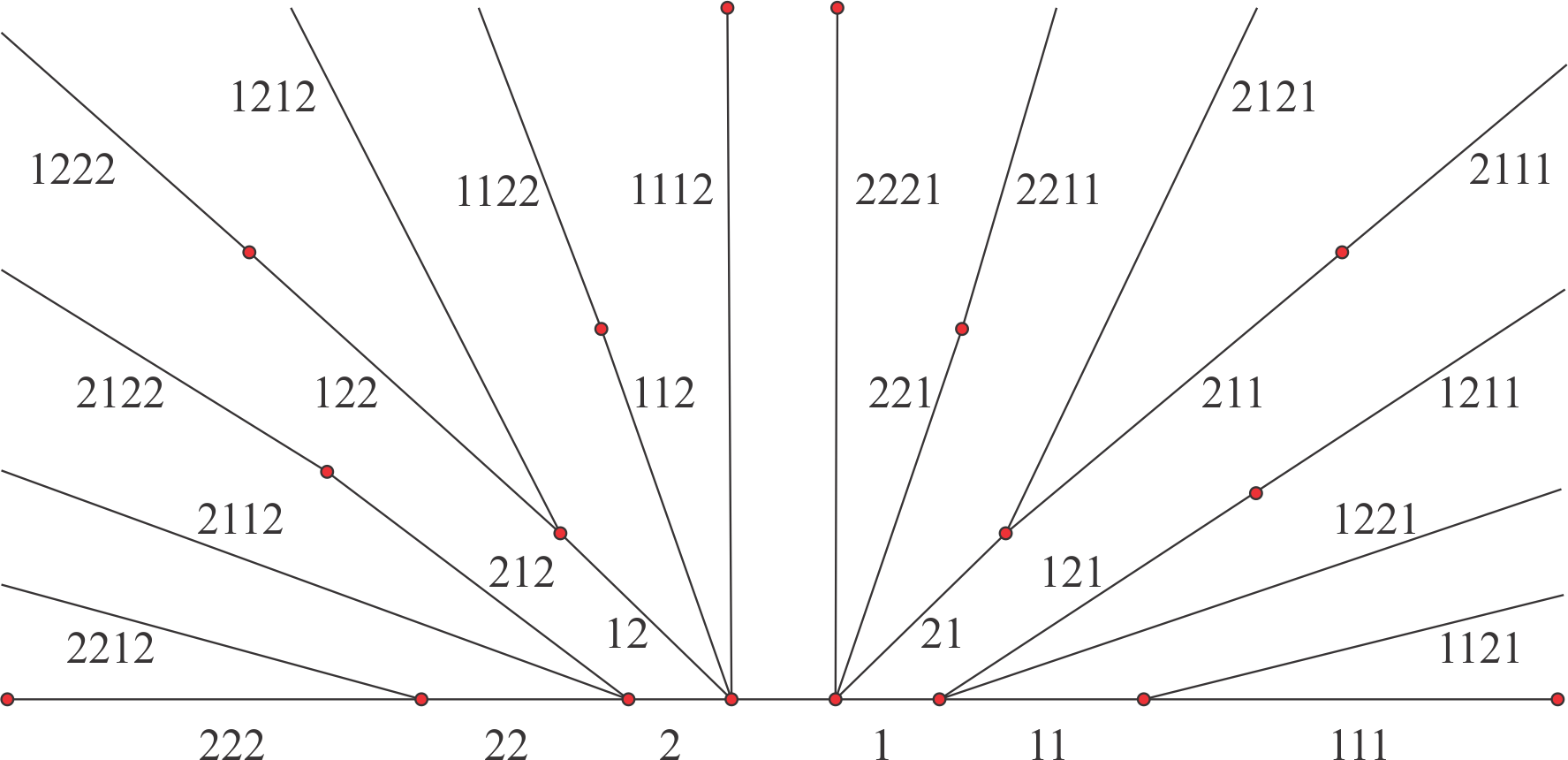}

\caption{A portion of the f-manifold of Example~\protect\ref{ex:binary}.
Each leaf $l(\theta)$ is labeled by $-\theta$. Note that
there are additional leaves, with longer labels, incident with some of the
branch points in the diagram, that are not included in the diagram.}
\label{fig:Example1}
\end{figure}

\begin{Example}
\label{ex:binary} The attractor of the IFS $\{ {\mathbb{R}} ; \, f_1,\,  f_2 \}$
where $f_1(x) = \frac12 x$ and $f_2(x) = \frac12 x + \frac12$, is the unit
interval $A = [0,1]$. The f-manifold is a connected, branched $1$-manifold.
The set ${\mathbb{K}} = l(\varnothing)$ is homeomorphic to~$A = [0,1]$.
According to Theorem~\ref{sevenpttwothm}, all other leaves are homeomorphic
to half open intervals. Each panicle, a union of leaves, is homeomorphic to
a closed interval. Routine calculation shows the following. The length of
the projection of each leaf, under~$\widehat \pi$, is a power of~$2$. Each
point $\alpha \in \mathbb{L}$ is contained in a neighborhood that is
homeomorphic to the open interval~$(0,1)$ or, if it is a~\textit{branch point}, to the union of a countably inf\/inite number of copies of~$[0,1)$ where
only~$\alpha$ is common to the copies. The branch points are isolated. More
exactly, the branch points are the points of $\mathbb{L}$ whose projection
under $\widehat \pi$ are of the form~$2^n$, $n\geq 0$, or $-2^n +1$, $n\geq 0$. Each sheet of $\mathbb{L}$ is mapped isometrically by $\widehat
\pi $ to either ${\mathbb{R}}$ or the interval $[0,\infty)$ or the interval $(-\infty,1]$. For example (recall that the ``bar'' means inf\/inite
repetition)  the sheet ${\mathbb{L}}_{\overline{-2 \, -1}}$ is homeomorphic
to ${\mathbb{R}}$, and ${\mathbb{L}}_{\overline{-2}}$ is homeomorphic to $(-\infty, 1]$. See Fig.~\ref{fig:Example1}.
\end{Example}

\begin{figure}[t]
\centering
\fbox{\includegraphics[width=3cm]{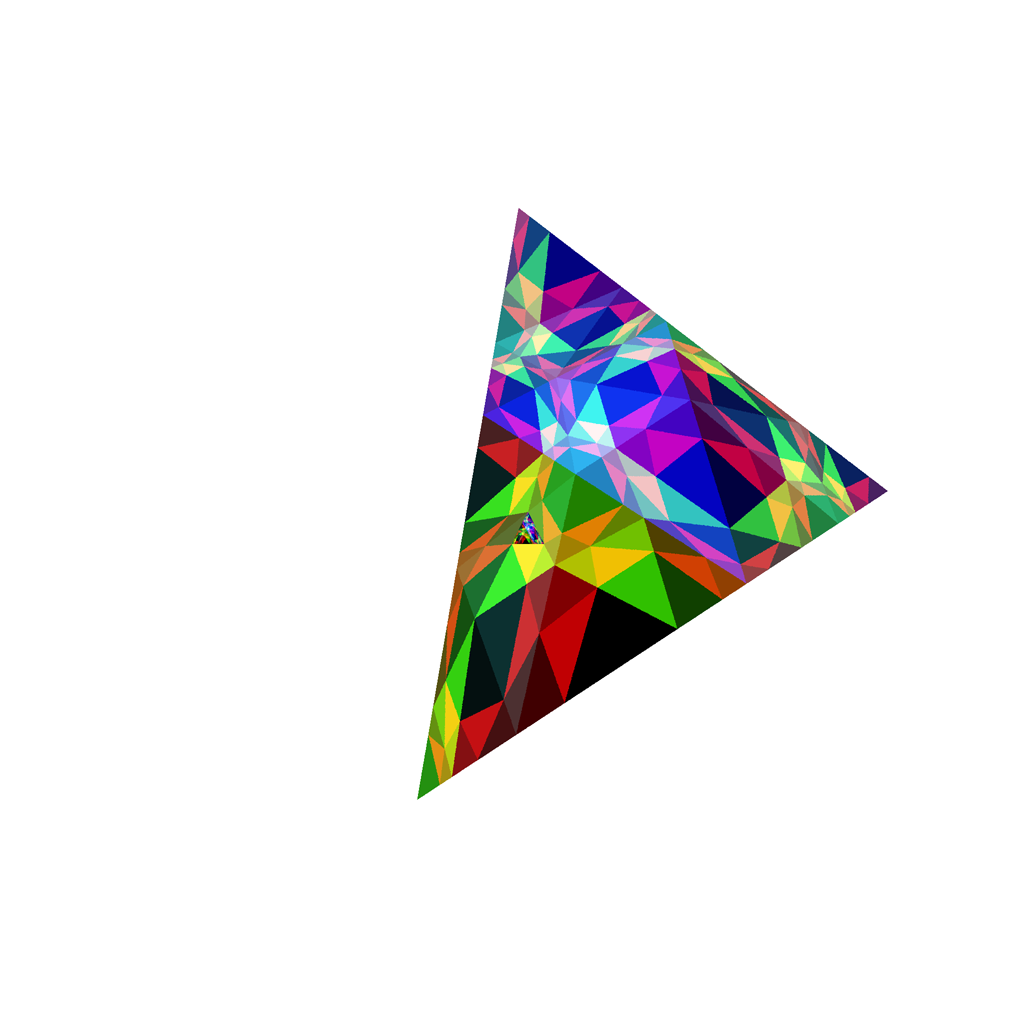}}
\quad \fbox{\includegraphics[width=3cm]{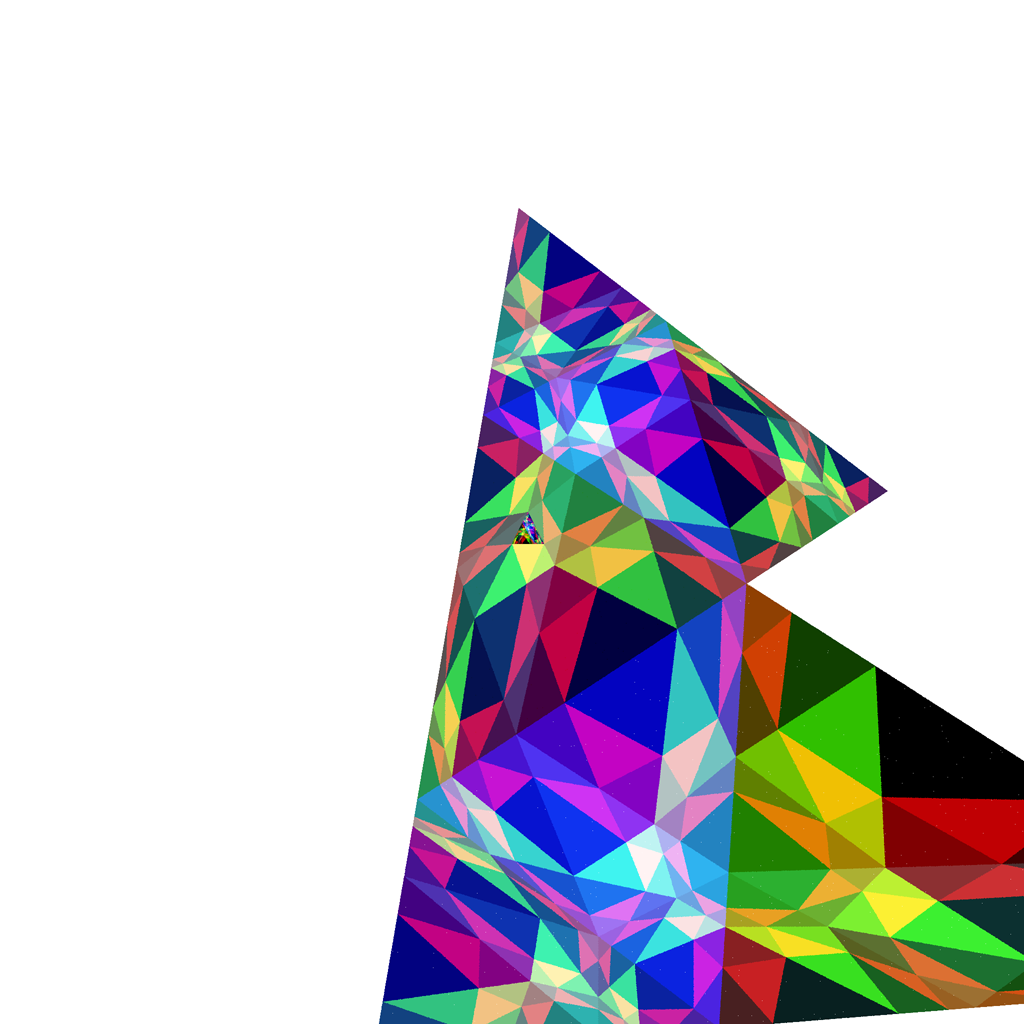}}

\fbox{\includegraphics[width=3cm]{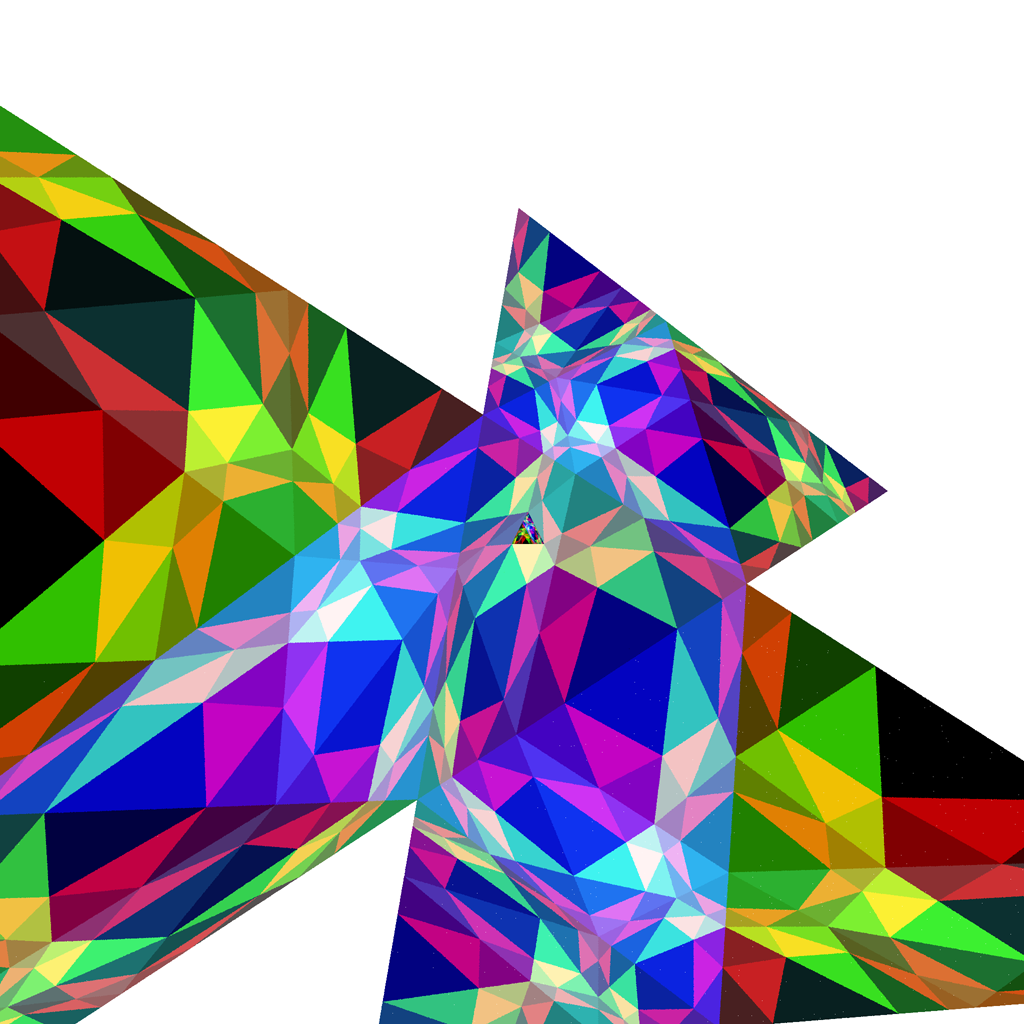}}
\quad \fbox{\includegraphics[width=3cm]{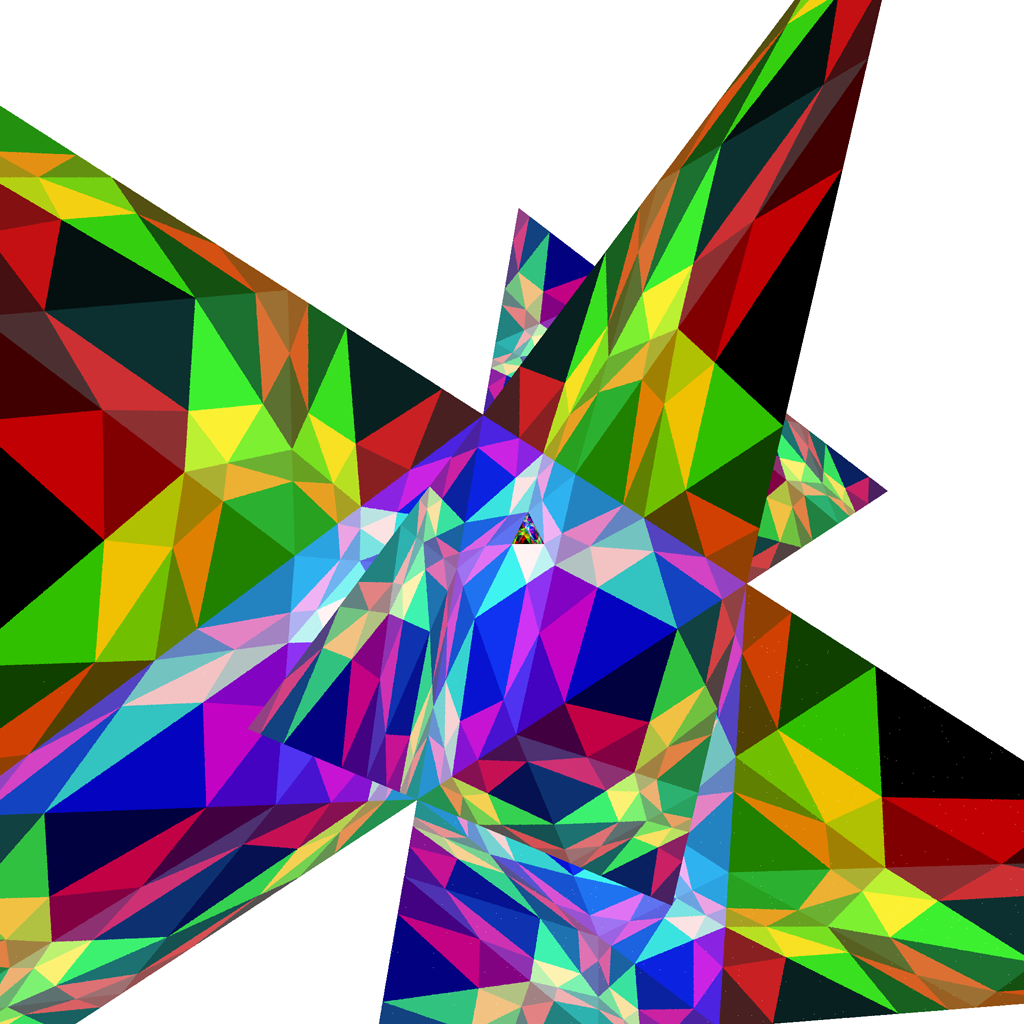}}
\caption{Fractal manifold under construction. See Example~\protect\ref{ex:manifold}.}
\label{fig:manifold008}
\end{figure}

\begin{Example}
\label{ex:manifold} While it is easy to visualize the branched fractal
manifold in the above example, it is harder to picture it in the following
example, illustrated in Fig.~\ref{fig:manifold008}. Here the IFS $\mathcal{F}=$ $\{\mathbb{R}^{2};\, f_{n}, \, n=1,2,3,4\}$ is af\/f\/ine, with the following
property. There is a triangle $ABC$ with points $c\in AB$, $a\in BC$, $b\in
CA$, where $XY$ is the line segment joining the points~$X$ and~$Y$, and the
triangles $abc$, $ABC$, $Abc$, $aBc$, $abC$ are non-degenerate. Moreover, $f_{1}(ABC)=Abc$, $f_{2}(ABC)=aBc$, $f_{3}(ABC)=abC$, and $f_{4}(ABC)=abc$.
IFSs of this kind are discussed in~\cite{monthly}. This IFS possesses a
unique attractor, the f\/illed triangle with vertices~$ABC$. This attractor is represented, in each of the four panels of Fig.~\ref{fig:manifold008}, by a very small multicoloured triangle located near the center of each of the
four pictures; directly below it is a small yellow triangle.  If $\theta\in\mathbb{I}_{+}$ is reversible, then $\widehat{B}_{\theta}$ $=\mathbb{R}^{2}$. The branched fractal manifold consists of non-denumerably many copies of $\mathbb{R}^{2}$ glued together appropriately. Each region of glue is
triangular.
The top left panel illustrates $B_{ijk1}$ for a f\/ixed choice of $ijk$. It
comprises 256 copies of $A$, namely the sets $f_{i}^{-1}\circ
f_{j}^{-1}\circ f_{k}^{-1}\circ f_{1}^{-1}\circ f_{m}\circ f_{n}\circ
f_{o}\circ f_{p}(A)$ for all $m,n,o,p\in\{1,2,3,4\}$, each copy in a colour
that is supposed to be dif\/ferent from its neighbours. The top right panel
illustrates $B_{ijk1}\cup B_{ijk2}$, where each of the subtriangles $f_{i}^{-1}\circ f_{j}^{-1}\circ f_{k}^{-1}\circ f_{2}^{-1}\circ f_{m}\circ
f_{n}\circ f_{o}\circ f_{p}(A)$ for all $m,n,o,p\in\{1,2,3,4\}$ is coloured
so that it is dif\/ferent from its neighbours. The two images may be thought
to be glued together on $B_{ijk}$. Similarly, the bottom left panel
illustrates $B_{ijk1}\cup B_{ijk2}\cup B_{ijk3}$ and the bottom right panel
illustrates $B_{ijk1}\cup B_{ijk2}\cup B_{ijk3}\cup B_{ijk4}$.
\end{Example}

\section{Related work}\label{sec:literature}

Fast basins, the extended code space, and branched fractal manifolds were f\/irst introduced in the earliest version of this paper~\cite{BLV}, which also concerns the notion of fractal continuation, f\/irst introduced in~\cite{continuations}.

The notion of fractal continuation is discussed in~\cite{BV8,continuations}, in the context of fractal interpolation, and in~\cite{tiling} in the context of tiling theory. The scope of the present work is much wider,
most particularly with regard to the def\/inition and topology of fast basins
and branched fractal manifolds. In~\cite{mfbkl}, a work in progress, conditions under which certain topological and geometrical properties that the fast basin inherits from the attractor of the IFS are established. These properties include the Hausdorf\/f and topological dimension, local connectivity and $\sigma$-porosity.

Notions related to but distinct from fast basins and branched fractal
manifolds occur  in the literature, mainly in the context of IFSs of
similitudes with the open set condition (OSC).
In particular, some fast basins are fractafolds as def\/ined by Stritchartz
\cite[Section~5.4]{strichartz1}, and some are  examples of fractafold
bundles as def\/ined by Ionescu and Kumjian~\cite{ionescu}. We discuss these
related notions in more detail in the following paragraphs.

In \cite{strichartz0} Strichartz def\/ines and establishes some properties of
fractal blowups (called fractal continuations here) for an IFS consisting of
similitudes that obeys the OSC. Fractal blowups for the classical Cantor
set, the Sierpinski triangle, and the Koch curve are investigated, and it is
noted that various dif\/ferent IFSs, with the Koch curve as the attractor,
yield dif\/ferent blowups. It is pointed out that in general a given IFS may
yield a nondenumerable set of non-homeomorphic blowups, where the relative
topology is assumed. Strichartz is focused on post critically f\/inite (PCF)
fractals because energy forms and Laplacians can be def\/ined on their
blowups. In \cite[Section~5.4]{strichartz1} Strichartz def\/ines a fractafold
to be a countable union of PCF fractals, glued together on their boundaries
in such a way that every point in the fractafold has a neighborhood that is
conjugate to a neighborhood of a point on the attractor; this notion
subsumes fractal blowups. In his context, Strichartz considers a fractafold
to be a generalization of the notion of a manifold, with the role of
Euclidean space (or half-space) being played by the attractor. In various
works, see for example~\cite{strichartz}, he has shown how energy forms and
Laplacians on PCF fractals can be extended to fractafolds.

Some branched fractal manifolds are related to fractafold bundles, def\/ined
by Ionescu and Kumjian~\cite{ionescu}. They consider IFSs of the form of $\mathcal{F}$ with $0<r_{i}\leq R_{i}<1$ and $r_{i}d(x,y)\leq
d(f_{i}(x)$, $f_{i}(y))\leq R_{i}d(x,y)$ for all $x,y\in \mathbb{X}$ and all $i=1,2,\dots ,N$, so that $\mathcal{F}$ is invertible and has a unique attractor~$A$. They def\/ine the associated fractafold bundle $L$ to be $\{\omega\times
B_{\omega}\colon\omega\in \mathbb{I}_{+}\}$, but with a certain inductive limit
topology in place of the relative topology. They use the projection of $L$
onto~$\mathbb{I}_{+}$ in order to build a Renault--Deaconu groupoid which
exploits the fact that if tails of two strings $\alpha,\beta\in \mathbb{I}_{+}$ are the same, and~$\mathcal{F}$ obeys the OSC, then~$B_{\alpha}$ and~$B_{\beta}$ are homeomorphic. This enables the establishment of properties of
an associated $C^{\ast}$-algebra on~$L$. The point of view in~\cite{ionescu}
is very dif\/ferent from that of the present work because it concerns the
projection of~$L$ onto $\mathbb{I}_{+}$ rather than into~$\mathbb{X}$.

\subsection*{Acknowledgements}

We thank Alan Carey for interesting discussions and helpful comments. We
thank Louisa Barnsley for technical help and suggestions. We thank Krystof Le\'sniak, who worked with us on the f\/irst version of this paper, namely \cite{BLV}, for many comments and suggestions. We thank the anonymous referees for helpful comments. This work was partially supported by a~grant from the Simons
Foundation (\#322515 to Andrew Vince). It was also partially supported by a grant from the Australian Research Council (\#DP130101738).

\pdfbookmark[1]{References}{ref}
\LastPageEnding

\end{document}